\NoBlackBoxes
\magnification =1100
\input epsf
\documentstyle {amsppt}
\def\a{\alpha}
\def\b{\beta}
\def\g{\gamma}
\def\G{\Gamma}
\def\d{\delta}
\def\D{\Delta}
\def\e{\epsilon}

\def\z{\zeta}

\def\i{\iota}

\def\l{\lambda}
\def\L{\Lambda}

\def\p{\pi}

\def\r{\rho}

\def\Si{\Sigma}

\def\Z{\raise.5ex\hbox{$\chi$}}

\def\bq{\Bbb Q}
\def\br{\Bbb R}
\def\bz{\Bbb Z}
\def\bc{\Bbb C}
\def\ca{\bc F}
\def\ra{\br F}
\def\rp{\br P(2)}
\def\cp{\bc P(2)}
\def\cf{\Cal F}
\def\rf{\br\cf}

\def\rF{\br F}
\def\cF{\bc F}
\def\rhf{\br \hat \Cal F}

\def\rhf{\br \hat \Cal F}

\def\rpt{\br P(3)}
\def\cy{\Cal Y}
\def\ct{\Cal T}
\def\cm{\Cal M}
\def\cP{\Cal P}
\def\cL{\bc L}
\def\J{\Cal J}
\def\ef{\eusm F}
\def\em{\eusm M}

\loadeusm
\topmatter
\title Floppy Curves \\ with Applications to Real Algebraic
Curves \endtitle

\author Patrick M. Gilmer \endauthor
\leftheadtext{Patrick M. Gilmer}%
 \address
 Department of Mathematics,
 Louisiana State University,
Baton Rouge, LA 70803 \endaddress
  \email gilmer\@ math.lsu.edu\endemail


 \date April 16, 1997\enddate

\thanks This research was supported by grants from
the Louisiana Education Quality Support Fund \endthanks

\abstract We show how one may sometimes perform singular ambient surgery on the
complex locus of a real algebraic curve and obtain what we call a floppy curve.
A floppy curve is a certain kind of singular surface in $\Bbb CP(2),$ more
general than the complex locus of a real nodal curve. We derive analogs for
floppy curves of known restrictions on real nodal curves. In particular we
derive analogs of Fiedler's congruence for certain nonsingular curves and
Viro's inequalities for nodal curves which generalize those of Arnold and
Petrovskii for nonsingular curves. We also obtain a determinant condition for
certain curves which are extremal with respect to some of these equalities. One
may prohibit certain schemes for real algebraic curves by prohibiting the
floppy curves which result from singular ambient surgery. In this way, we give
a new proof of Shustin's prohibition of the scheme $1<2\coprod 1<18>>$ for a
real algebraic curve of degree eight.
\endabstract

\endtopmatter

\document
\head Introduction \endhead

A real algebraic curve $\ra$ in $\rp$ of degree $m$ is the
set of zeros
of a real homogeneous polynomial $F$ of degree $m$ in three
variables. We
will refer to $\ra$ as a real curve. It is nonsingular if the gradient
of $F$ is
nonzero on $\ra$. In this case $\ra$ is a collection of smoothly
embedded
simple closed curves in $\rp$.
We let $\ca$ be the set of zeros of $F$ in ${\Bbb C}P(2)$, and refer
to  $\ca$ as
a complex curve.  If
$\ca $ has only nodal (ordinary double point) singularities and
no multiple components,  then we call $\ra$ a nodal curve, and $\ca$
a complex
nodal curve.  ${\Bbb R}A$
will consist of a collection of immersed curves in the plane with
only double point singularities which intersect transversely away
from these singularities together with a finite set of points.   The
possible
ambient isotopy types of such subsets of
${\Bbb R}P(2)$ is interesting. In 1978, Viro \cite{V3} found the
analogs for even
degree nodal real algebraic curves of the strengthened Petrovski\u i
and the strengthened  Arnold inequalities  \cite{W,
(7.4)},  \cite{V2,
(3.4), (3.5), (3.6), (3.7)}. The original inequalities applied to
nonsingular curves
of even  degree. (At the same conference Kharlamov and Viro
independently stated a version of the original Petrovski\u i inequality for
more
general singular curves \cite{Kh1} \cite{V3}).  Kharlamov and Viro have also
studied
analogs of the Gudkov-Rokhlin
congruence for real nodal curves \cite{KV}(they also consider more
general
singularities). Fiedler has also given some generalization of his
congruences
for nonsingular curves \cite{F1} which also applies to  nodal curves
\cite{F2}.
\par
Viro introduced the notion of a flexible curve \cite{V2,p.59}. This is
a surface in
$\cp$ which is homologous to $m [\bc P(1)]$ satisfying three further
conditions
which are automatically satisfied by a complex curve. Viro observed
that many
known restrictions apply equally well to flexible curves.
We will define a  similar weakening of the notion
of a  nodal curve.  As our notion allows considerably  more
freedom than one might expect, we will refer to this new kind of
curve as a floppy curve.  In this paper
we show that the  existence of a  real nonsingular curve frequently
implies the existence of related  floppy curves.  If these
floppy curves can be prohibited, so are the
original real nonsingular curves.
\par
We further generalize Viro's  inequalities for nodal curves  so that they apply
 to floppy
curves as well.    Our results
\cite{G1}
relating nonsingular curves to link cobordisms also apply to floppy
curves. Thus the restrictions on cobordisms given using signatures
and Arf
invariants given in \cite{G2} also apply. Thus many theorems are
available to
prohibit floppy curves. In particular we derive Theorem (3.3) which gives
restrictions on floppy curves based on the Arf invariants of associated links
in the tangent circle bundle of $RP(2).$  Theorem (3.3) when applied to
nonsingular curves yields Fieder's congruence \cite{F1}. Fiedler observes that
his result \cite{F2} also applies to
flexible nodal curves. Theorem (3.3) also probably implies
\cite{F2}, but we have not checked this yet. However Theorem (3.3) applies to
some curves to which \cite{F2} does not apply. To illustrate our methods, we
use Theorem (3.3) to give a new proof of Shustin's  prohibition of the scheme
$1<2\coprod 1<18>>$ for a real algebraic curve of degree eight.

\par
Our method of constructing floppy curves makes use
of Fiedler's concept of a chain of ovals \cite{F3}, and Viro's idea
\cite{V2,p67},\cite{K-S,\S 5} of using chains of ovals to find
membranes for
$\ca$.  Our method requires that one deduce from the isotopy class
of a
hypothetical real non-singular curve some set of chains of ovals
associated to a
given pencil of lines (or a list of possibilities for such sets of
chains of ovals, one
of which must occur).
There are many arguments in the literature \cite{V4} \cite{K-S},
\cite{K1},
\cite{K2} for deducing chains of ovals.

There are other reasons to study nodal curves. The union of a line and
say a
nonsingular curve of degree $m'$ is a nodal curve of degree $m'+1$ .
As the
complement of a line is affine space,
the study of such arrangements is the study of affine curves of
degree $m'$. This question has received a lot of attention
\cite{F3,Theorem 2},\cite{P},\cite{KS}.  Besides its intrinsic
interest, one reason to study this case is that it is related to
smoothings of non-degenerate singular points \cite{S1},\cite{KS}.
This in turn is important for the construction of nonsingular curves
by Viro's method \cite{V5}.
Another way the above theorems may prove useful is as follows.
Sometimes one can use B\'ezout's theorem to show an isotopy class
of
a nonsingular curve implies an isotopy class for the arrangement (or
one of several possible arrangements) made up of a given curve and
some auxiliary lines (or conics). Then one can rule out the original
curve by showing the resulting arrangements are impossible
\cite{V4}, \cite{K1}, \cite{K2}.  We provide new methods
which may be added to the arsenal of methods which can be used to
rule out the resulting arrangements.

Section one contains preliminaries on nodal curves and the
definition of a floppy curve.
Section two describes how floppy curves can be constructed by ambient singular
surgery on real algebraic curves.
In section three, we discuss link cobordisms related to flexible nodal curves.
 We also give our generalization of Fiedler's congruence.  In section
four, we state the generalized
Viro inequalities for floppy curves.
Section five contains a further restriction on floppy curves involving the
determinant of some matrices used in the statement of the generalized
Viro inequalities for floppy curves.
Sections six through nine provides a complete proof
the generalized
Viro inequalities for floppy curves, as well as the further determinant
conditions. This should prove useful even for those only interested in Viro's
original theorem as \cite{V3} is an announcement and contains no proofs.
 As we were completing an earlier version of this paper in April
1996, we learned of a recent paper by Finashin \cite{Fi} which derives the
Viro-inequalities by somewhat different approach.   We thank Viro for useful
conversations.
The basic method of proof which we use is the same as Viro's which
in turn is based on Arnold's original method \cite{A}.
We were greatly influenced by Wilson's exposition and extensions \cite{W}.

We let $\beta_k(X,\Bbb Z_p)$ denotes $\dim H_k(X,\Bbb Z_p)$.
Let $T_p(M)$ denote the real tangent space of a manifold M at a point $p
\in
M$.  Let $T$ denote the involution on $\cp$ given by complex
conjugation,
i.e. $T[z_0:z_1:z_2] = [\bar z_0:\bar z_1:\bar z_2]$.

\head \S 1 Nodal Curves and Floppy Curves \endhead

Let $G$ be a complex polynomial factor of $F$, and let
$\bc G$ and
$\br G$ denote the associated zero sets in $\cp$ and $\rp$.  If
$G$ is irreducible over $\bc$,  we call $\bc G$ an
$\bc$-irreducible component of $\ca$.  If $G$ is real, and is irreducible
over $\br$, $\bc G$ is called an  $\br$-irreducible component.  If $\bc G$ is
$\br$-irreducible but not $\bc$-irreducible,  then $G= g\bar g$ where  $g$ is a
complex polynomial irreducible over $\bc$ and
$\br G$ must consist of isolated points which are nodal
singularities for $\bc G$. Otherwise $\bc g$ and $\bc \bar g$ would be distinct
irreducible curves with an infinite number of intersections.

 \par
Let $m$ be a non-negative integer.   By a floppy curve of
degree $m$, we mean the image $\cf$ of an immersion $j$ of a closed
oriented surface $\hat \cf$ in $\cp$  with only ordinary double point
singularities which is invariant under $T$. Note that $T_{|F}$ lifts to an
involution on
$\hat \cf,$ which we also denote by $T.$ Let $\rf$ denote
$\cf \cap \rp$ and $\rhf$ the inverse image of $\rf$ in $\hat \cf$.    $\rhf$
is a disjoint collection of simple closed curves and some isolated points.
$\rf$
consists of a collection of immersed circles in $\rp$ together
with  a finite number of isolated double points. Each of these immersed circles
is
called a constituent of $\rf$. Let $\rf^*$  denote $\rf$ with any
isolated points deleted, and $\rhf^*$ its inverse image in $\hat \cf$ . We
insist that the following  conditions hold.
 \roster
\item[1] $\cf$ represents the class $m[\bc P(1)] \in H_2(\cp).$
\item[2] For each $p\in \rhf^*$, there is a line $L_p$ in
$ T_{j(p)}(\rp)
\subset T_{j(p)}(\cp)$ such that
$j_*(T_p(\rhf^*)) \oplus i(L_p) = j_*(T_p(\hat \cf)$
\item[3] If $p \in \rhf -\rhf^*$, we insist that $j_*(T(\hat \cf)_{|p})$ is the
realification
of  a complex line in the complex tangent bundle of $\cp$  and the orientation
on
$\cf$ agrees with the orientation arising from the complex structure.
\endroster
Note that for each $p\in \rhf^*$, there can be at most one  line $L_p$
satisfying
the above condition.
We denote the function $p \rightarrow  L_p$ by $\Cal L_{\rf}$.
If $j(p_1)=j(p_2)$ is a double point of $\rf^*$, then $L_{p_1} \ne
L_{p_2}$. The orientation on
$\hat \cf$ induces an orientation $\eufm o_{\rf}$ on $j_*(T^{\br}_p(\rhf^*)
\oplus
i(\Cal L_{\rf}(p))$.

\par We write $\rf$ when we mean  the triple $(\rf,\Cal L_{\rf},\eufm o_{\rf})$
and
call $\rf$ the real part of $\cf$.  We refer to $\rf$ as a real floppy curve.
Thus our
real floppy curves come equipped with a line field $\Cal L_{\rf}$, and an
orientation
$\eufm o_{\rf}$ as above. Note that condition (3) above forces all isolated
real double
points to be positive double points of $\cf$.
\par
Sketching a line field can be quite tedious. However there is an easy way to
encode in a drawing all  the data of a triple $(\rf,\Cal L_{\rf},\eufm
o_{\rf})$, up to isotopy.
Note that we can isotope the line field so that it is tangent to $\rf^*$ in a
neighborhood of every double point. Then we may further isotope the line field
so
that it is tangent
everywhere except along some short segments where it looks like the field
illustrated in Figure 1a.   We call these segments flops. It is not hard to see
that one may isotope the floppy curve $\cf$ so that
$\rf$ has this new field. In this situation we say that $\Cal L_{\rf}$ is
mostly tangent to $\rf$. We can indicate in a
drawing of $\rf$ that the line field has a flop by replacing such a
segment by
a segment with a ``cusp'' as illustrated in Figure 1b.
\midinsert
$$\vbox{\epsffile{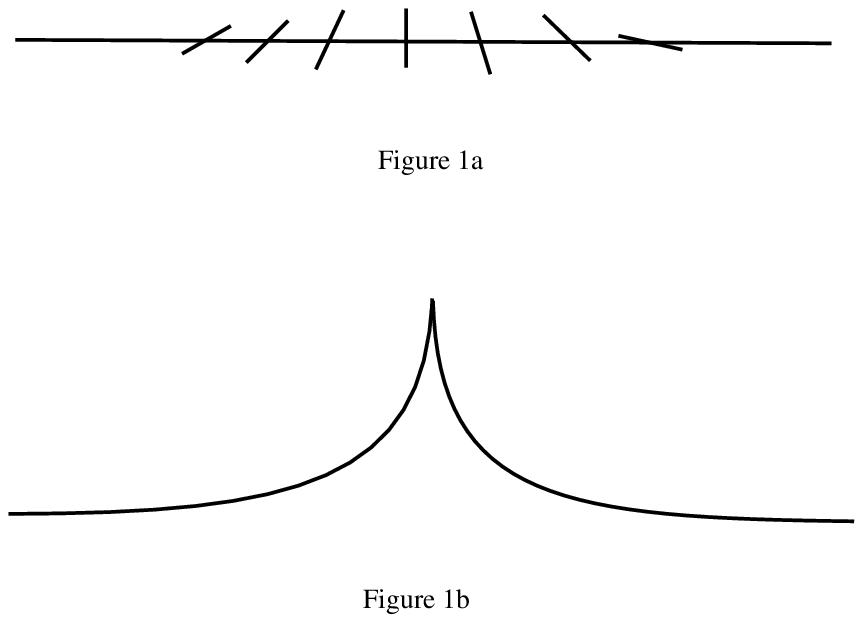}}$$
\endinsert
 So as not to encourage any further associations with cusp singularities beyond
the visual associations, we will refer to these cusp-like segments of the
encoded drawing as flops. Notice that the direction of the
flop is determined by the line field.  At a point $p$ where the line field is
tangent to the
appropriate branch of $\rf$, the orientation given by a tangent vector  $v$
followed by $i v$ will either agree or disagree with $\eufm o$. This is
independent of our choice of $v$.  We say $\eufm o (p)$ is $\pm 1$, depending
on
whether these orientations agree or disagree. This function changes exactly
when one passes through a flop. It follows that  each
constituent of $\rf$ contains  an even number of flops. We can encode this by
drawing the segments where $\eufm o$ is positive thicker.  Thus we can encode
the
isotopy type of a real floppy curve with a collection closed curves in $\rp$
which meet transversely and each of which has an even number of flops.
Moreover the thickness of the curves should change at each flop as we travel
along each constituent of $\rf$ .
\par

A nonisolated real double point of $\rf$ will be a positive double point of
$\cf$ if and only if the value of $\eufm o$ agrees on the two branches which
cross at the double point. In our encoded pictures this mean both branches are
thick or both branches are thin.  If a floppy curve has a nonisolated real
double point, we may obtain
another floppy curve by resolving the double point in either
of the two possible ways as in Figure 2. Figure 2a shows the case where both
curves are thick, there is a similar case where both curves are thin.
 \midinsert
$$\vbox{\epsffile{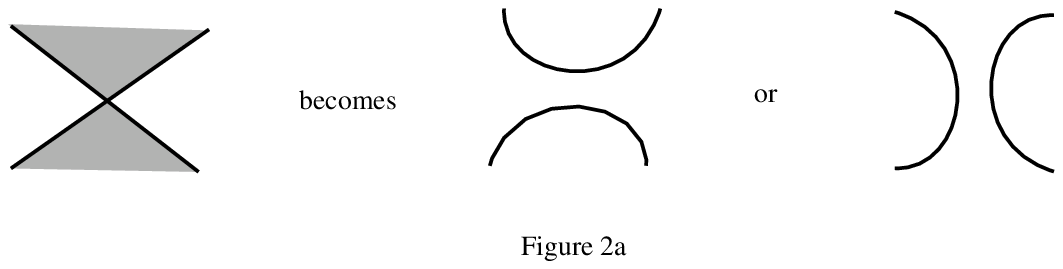}}$$
\endinsert
 \midinsert
$$\vbox{\epsffile{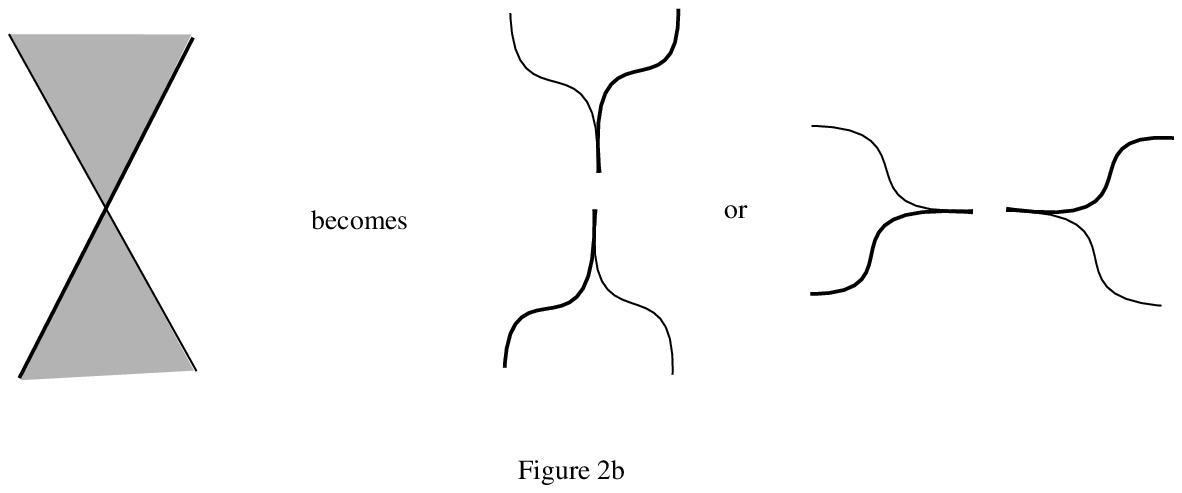}}$$
\endinsert

The topology of $\cf$ is modified under these moves by deleting a cone on the
disjoint union of two circles and gluing in an annulus. To describe this
precisely recall there is a orientation reversing diffeomorphism $\Cal D$ from
the tangent disk bundle of $\rp$ to its tubular neighborhood in $\cp$ given by
multiplication by $i$. One defines on the shaded regions in Figures 2a and 2b.
unit line fields.  In Figure 2a it is $\pm <x,-y>$, after we have chosen
coordinates which make the lines the $x$ and $y$ axes. In Figure 2b it is $\pm
<x,y>$. One completes these by an arc of lines located at the double points.
Under $\Cal D$, these describe a pair bands $I \times I$ in the sphere bundle
of $\rp$ which meet  $\cf$ intersected with the before pictures of along
$\partial I \times I$ and meets the local picture of the floppy curve given by
the first ``after'' picture of the ``after'' picture along $I \times \partial
I$. Each band is given by a unit vector field parallel to the line field
together with an arc of unit vectors located at the double point. On may then
replace $\cf$ intersected the disk bundle of the before picture by the two
bands union the two disks which are given by the  tangent unit disk bundles of
the first after picture under $\eufm D$.  To obtain the second ``after''
pictures one just considers the bands given by the above line fields defined on
the unshaded regions. These is are local moves and
do not change $\cf$ away from a small neighborhood of the point in question.
The next section will show that under certain circumstances the reverse
process of Figure 2a can be applied to a nonsingular real algebraic curves.
\par
{}From a floppy curve with some number of pairs of complex
conjugate double points, one may obtain a floppy curve with one fewer
pair of complex conjugate double points by
equivariantly resolving
the pair of double points.

Let $n(\cf)$  denote the number of connected components of $\hat
\cf$.   We will call the image of a component of $\hat \cf$ an
$\bc$-irreducible component of $\cf.$ The degree of a $\bc$-irreducible
component of $\cf$ is the nonzero multiple of the homology class of $\bc (P(1)$
that it represents in $H_2(\bc P(2)).$
Note that a $\bc$-irreducible component of $\cf$ is not necessarily a floppy
curve as it need not be invariant under $T.$
A floppy curve $\cf$ is said to be
$\br$-irreducible if either it is $\bc$-irreducible or $\cf$ consists of pair
of $\bc$-irreducible components which are interchanged by $T.$
If $\cf$ is $\br$-irreducible but not $\bc$-irreducible, then $\rf$ consists of
a  collection of isolated points.
$\cf$ is the union of a collection $\br$-irreducible curves, called the
$\br$-irreducible components.
Let $c(\cf)$ denote the number of $\br$-irreducible components of $\hat
\cf$. This is also the number of components of the orbit space $\hat \cf$ under
$T.$
If  a $\br$-irreducible component of $\cf$ is $\bc$-irreducible, we will call
it a strongly-irreducible component. If  a $\br$-irreducible component of $\cf$
is not irreducible, we will call it a weakly-irreducible component.

An strongly-irreducible component is called dividing if it becomes
disconnected when its real part and its non-real double points are deleted.
For such a curve, a choice
of one half of $\cf -\rf$ induces a semi-orientation on
$\rf^*$ called the complex orientation. We can indicate a complex orientation
on a real floppy curve by drawing  arrows  on each constituent of
$\rhf^*$ in a consistent way. $\cf$ is  called dividing, if each
strongly-irreducible component is dividing.

\par
If a floppy curve has an isolated real point  then one may at
will either delete it or replace it with a small empty oval with no flops and
$\eufm o$ identically one on the added oval. Again the topology of $\cf$ is
changed in this process by replacing  a cone on two circles with an annulus.
This would not be possible if we did not insist on the orientation part of
condition (4). Gudkov showed that  modifications similar to this one and Figure
2a could be
made to (non-flexible) nodal curves. See \cite {Gu, p12.}

If the
original curve is dividing and we resolve it by creating a small empty oval
then so is the resulting curve with an
induced complex conjugation which remains constant on the
unaffected
portions of $\rf$.  Suppose we have chosen a component of $\cf_i-\rf_i$ for
each strongly-irreducible component $\cf_i$ and also chosen a component of
$\cf_j$ for each weakly-irreducible component $\cf_j$.  Then the resolved curve
we obtain by replacing an
isolated double point of a dividing floppy curve with an empty oval comes with
such choices made in a natural way.   Thus we can speak
of a complex orientation for a dividing floppy curve which includes this
information. This can be indicated in a picture of $\rf$ by an arrow curved
around the isolated point showing the way an added oval would be oriented.
In the moves described by Figures (2) and (2b), if the original floppy curve is
dividing , then the resulting curve will
be if we make the above resolution compatible with the complex
orientation. Otherwise the resulting curve will not be dividing.

\proclaim {Lemma (1.1)}  $\cf$ is a floppy curve of even
degree if and only if $\rf^*$ is null-homologous modulo two in $\rp$.
\endproclaim
\demo{Proof} We may isotope $\cf$ slightly through floppy curves so that it
meets $\Bbb CP(1)$ transversely in points. The number of intersections is even
and all the nonreal intersections arise in pairs. The intersection of
$\rf^*$ with $\Bbb RP(1)$ is even.\qed
\enddemo

We emphasize that  $\cf$ may have both positive and negative double
points.   We let $d_{\pm}(\cf)$ denote the number of positive
(respectively negative)
double points and let $d(\cf)$ denote the total number of
double points.   Define

$$ \D(\cf) = \chi (\cf) -d_+(\cf)  +d_-(\cf) +m^2 -3m.\leqno(*)$$
Note that a nodal curve $\cF$ of degree $m$ is a floppy curve
of degree m.  Although a nodal curve will automatically satisfy the
following additional conditions which can be phrased topologically,
we do not require that floppy curves satisfy them :
\roster
\item[4] All the double points are positive, i.e. $d_-(\cF) =0$.
\item[5] The oriented tangent planes to $\cF$ at points in $ \cF\cap
\rp$ are complex lines in the tangent space of $\cp$ and the
orientation induced by the complex structure agrees with the given
orientation.
\item[6] $\D(\cF)=0$.
\endroster
To see that condition seven is satisfied. Perturb $F$ slightly
obtaining  $\check F$ with $\bc \check F$ is nonsingular. One has
$\chi
(\bc \check F) = 3m-m^2$, and  $\chi (\bc \check F)=\chi (\cF)-
d_+(\cf)$.  A floppy curve $\cf$ may be isotoped so that Condition (5) holds if and only if  $\rf$ has no flops
and $\eufm
o$ is identically $+1$.  It is because we do not require floppy curves to
satisfy conditions (4),
(5), and (6) that we did not use the name flexible to describe them. We define
a
flexible nodal curve to be a floppy curve which satisfies these last three
conditions. The flexible curves of Viro are then flexible nodal curves without
any double points.

We also
have
the following generalization of Harnack's theorem, and a related congruence for
dividing curves \cite{V2,(3.12)}. Let $r(\cf)$ denote the number double points
which lie in $\rp$ and are not isolated in $\rf$. So $d(\cf)= 2\nu(\cf)+
i(\cf)+r(\cf),$ where $\nu(\cf)$
denotes the number of pairs of complex conjugate double points
of $\cf$ in $\cp - \rp$. We will also let $\nu_+(\cf)$ ($\nu_-(\cf)$)
denotes the number such  pairs of positive (negative) double points.
Let $\ell (\cf)$ denote the number
of constituents of $\rhf^*$. Let $$h(\cf)=  4 c(\cf)-\chi(\cf) - 2 \ell (\cf) -
d(\cf).$$

\proclaim {Proposition (1.2)}
$h(\cf) \ge 0.$
This is an equality modulo two. For dividing curves it is an equality
modulo 4. If $h(\cf) = 0,$   then the curve is dividing. \endproclaim
\demo{Proof} Let $\tilde G$ denote $\cf$ with $\eufm n (\rf),$ the interior a
small closed
neighborhood of $\rf$ (invariant under complex conjugation), deleted. 
$\tilde G \cap
\overline{\eufm n (\rf)}$ consists of $2 \ell (\cf)+2 \i(\cf)$ circles and so
has zero
Euler
characteristic.  Thus $\chi(\tilde G) = \chi( \cf)-\chi(\rf)$.  Let $G$ denote
the orbit
space of $\tilde G$ under complex conjugation.
We have $\chi(\tilde G) =2\chi(\
G)$. $G$ is a surface with $\ell (\cf)+ \i(\cf)$ boundary components  with $\nu
(\cf)$
interior points identified. Let $\bar G$ denote the closed  surface with
$c(\cf)$ components that
obtained by capping off each boundary component with a disk and unidentifying
any interior points.
$$\chi (\bar G) = \chi (G) + \ell(\cf) + \i(\cf) +\nu (\cf).$$
The above inequalities simply say that $\chi(\bar G)$ must be an integer less
than
or equal to $2 c(\cf)$ and must be even if $G$ is orientable.  If $h(\cf) = 0,$
, then $\bar G$ is a collection of 2-spheres, and the curve must be dividing.
\qed\enddemo

 \head \S 2 Pencils of
Lines, and Floppy  Curves \endhead
Fiedler introduced pencils of lines to the study of real algebraic
curves \cite{F3}.   Let
$\rF$ be a nodal real algebraic curve. Let $P \in \rp$ be  point which
does not lie on $\rF$ or
on any  inflectional tangent lines to  $\cF$ or on any real
line which is tangent to $\cF  -  \rF$ or on any real line
which  passes through a  node on $\cF  -  \rF$. Here a real
line is simply a degree one curve with real equation.  The set of such
$P$
is open and dense as the set of illegal points is a proper real
algebraic set. Let $\Cal L (P) = \{
 \br L_u\}_{u \in \br P(1)}$  be the pencil of all real lines through
a point $P$. Here $L_{[u_0:u_1]}= u_0 L_0 + u_1 L_1$ where $\br
L_0$ and $\br L_1$ are distinct real lines through $P$.  For $t \in \br$,
we let $L_t$  denote $L_{[1-t:t]}$. This is consistent with the above
notations for $L_0$ and $L_1$. We call $\{L_t | t \in [0,1]\}$ a segment of the
pencil of lines joining $L_0$ to $L_1$. There are only two such segments, one
of which is picked out by a choice of sign for the defining equations of $L_0$
and $L_1$.
Choose a  complex orientation of each $\br L_t$ in a consistent way ``so that
the orientations pass into one another under the natural isotopy''  and denote
by $\bc L{_t}^+$ the closure of the component of $\bc L_t - \br L_t$ which
induces the complex orientation. The index of a point of tangency of $\rF$ with
a line in $L \in \cP$ is the Morse index of the function which sends $x$ to $t$
whenever $x \in L_t$.  Consider a pair of tangencies of lines in $\cP$ with a
component $C$  of $\rF$ where one
has index one and the other has index zero. We say that such a pair of
tangencies is inessential if and only if the  complex orientation on the
tangent lines induce the same orientation on $C$. Otherwise  such a pair of
tangencies is called essential. See \cite{V4,p.413}.
Following Viro
\cite{V2} we denote  the closure of $\cF \cap((\cup_{L \in \Cal L (P)
}\bc L)  -  \rF )$ by $S_P(F)$.
$S_P(F)$ is a smooth closed one dimensional submanifold of $\cp$
lying on $\cF$.
$S_P(F) \cap \rF$ consists of the points of tangency of $\Cal L (P)$ with
$\rF.$
\par
In \cite{V2,p.67}, Viro announced some new prohibitions for curves of degree
eight obtained by a new method of using embedded membranes  for $\cF$ with
boundary the union of some of the components of $S_P(F)$.  By an embedded
pseudo-membrane (e-p-membrane for short) for $\cF$ we mean an embedded surface
$\cm$ such that $\partial \cm \subset \cm \cap \cF$. By an embedded membrane
(e-membrane for short), we mean an e-p-membrane whose interior is transverse to
$\cF$.   Viro's method was first described in detail in \cite{KS,\S 5} where it
is applied
to the study of affine curves of degree six. An e-membrane can be found in the
following situation which is  more general than that described in  \cite{KS}.
Still we follow the exposition in \cite{KS} quite closely.  There is a misprint
in the English translation \cite{KS,p509}:$\bold C L \setminus \bold R L$
should read $\bold C L \cap \bold R L$. Note in \cite{KS}, the membrane is used
to construct homology classes in a branched double cover along $\cF.$ We plan
to use it to perform surgery on $\cF.$ This accounts for some differences in
our treatment.
  Let $\cP$ be a segment of the pencil of lines joining two lines through $P$:
$L_0$ and $L_1$.
  Let $\cP^+ = \cup_t \bc L{_t}^+$.
and $\br \cP = \cP^+ \cap \rp$. Note that $\cP^+$ is homeomorphic to $D^2
\times I / \{p\} \times I$ where $p$ is a point on the boundary of $D^2$.   We
assume that the following conditions
hold:
\roster
\item $\cP$ contains an even number of simple tangents to the curve $\cF$.
These tangents correspond to tangents of $\br L_t$ and $\br F$ in $\rp$.  The
remaining lines in $\cP$ intersect $\cF$ transversely in $2k$ points, of which
at least $2k-4$ are real.
\item $L_0$ and $L_1$ are tangent to $\cF$ and intersect it only in real
points.
\item The points of tangencies of the lines in $\cP$ with $\rF$ can be grouped
into pairs with each pair joined by an arc,  such that the arcs are disjoint
and every line $L_t$ in
$\cP$ which is not a tangent to $\rF$ transversely intersects exactly $2-s/2$
arcs, where  $s$ is the number of imaginary points in $\bc L_t \cap \cF$. Note
$s$ is either zero, two or four.  The
arcs must be disjoint and  must be one of following two types. The first type,
called inessential, join an
inessential pair of tangencies on the same oval, meets $\rF$ only at its
endpoints, is tangent to $\rF$ at its endpoints  and is a small perturbation of
an arc on this oval joining these points.  The second kind, called essential,
is
transverse to $\rF$, misses the double points of $\rF$ and either joins points
on  distinct
ovals or a pair of essential tangencies.  \endroster
\par
Let $\{ \L_i\}$ denote the set of inessential arcs and let
$\{ \J_j\}$ denote the set of essential arcs.  Let $S= (\cup \L_i) \cup (\cup
\J_j) \cup (\cP^+
\cap (\cF-\rF)$.  Each disk $\cL_t^+$ for $t\ne 0, 1$
 contains exactly two points of $S$ and $\cL_0^+$ and $\cL_1^+$ contains
exactly one point of $S$. Thus $S$ is a circle, which may be spanned by a
disk $\D$ whose interior consists of one arc in each $\cL_t^+$ for $t\ne 0, 1$.
Note $\Cal M =
\D \cup T(\D)$ is a planar surface with Euler characteristic
two minus the number of arcs, and is an ep-membrane for $\cF$ with boundary
some
of the components of $S_P(F)$.
\par
 In deducing the occurrence of the above set up, from the topology of a
hypothetical isotopy type for $\rF$, one usually has no control over how many
inessential arcs are needed. However one can find a nearby e-membrane  $\Cal
M'$ with Euler characteristic two minus the number of essential arcs. First we
consider the e-p-membrane $\Cal M''$ constructed as above only replacing
the inessential arcs by the unperturbed arcs $\G_i$ on the components joining
inessential pairs. Let $\G = \cup \G_i$, and let $S'_P(F)$ be $S_P(F)$ surgered
in $\cF$ along the $\G$. We will construct $\Cal M'$ to have boundary the union
of some of the
components of  $S'_P(F)$. We may isotope (see below) $\Cal M$ to say $\Cal M''$
 in  a neighborhood of  $\G$ so that $\Cal M'' \cap \cF$ is a tubular
neighborhood of
 $\G$ in $\cF$.  Let $\Cal M'$ denote
$\Cal M'' -\text{Int}(\Cal M'' \cap \cF)$ smoothed.
\par
We must describe the isotopy of  $\Cal M$ to $\Cal M''$ precisely. Recall that
multiplication by $i$ defines an orientation reversing diffeomorphism from a
tubular neighborhood $\Cal N$ of $\rp \subset \cp$ to the tangent disk bundle
of $\rp$ which sends $\cF \cap \Cal N$ to the total space of the line
bundle given by $\Cal L_{\rF}$. Moreover it sends $\Cal M\cap \Cal N$ to the
total space of the line bundle over $\cup (\L_i)\cup (\cup(\J_j))$ given by the
line field  which sends each point $p$ to the line joining $p$ to $P$ which we
denote by $\Cal L_P$. Our isotopy takes place completely inside of $\Cal N$.
Under
our diffeomorphism it completely lies over the disks $D_i$ which span the union
of  inessential $\L_i$ and the associated $\G_i$. We pick a line segment field
of fixed length over each $D_i$ which restricted to $\G_i$ is tangent to $\rF$
i.e. given by $\Cal L_{\rF}$  and which restricted to $\L_i$ is parallel to
$\Cal L_P$. See Figure 3.  Then the total space of this line field is a 4-disk
we push $M$ across to get $M''$.
It is important to note that $\Cal M$ lies completely in $\cP^+ \cup T(\cP^+)$,
and that $\Cal M'$ lies completely in $(\cP^+) \cup T(\cP^+)$ except for a
portion which lies over $\cup D_i$ under under the diffeomorphism.
\midinsert
$$\vbox{\epsffile{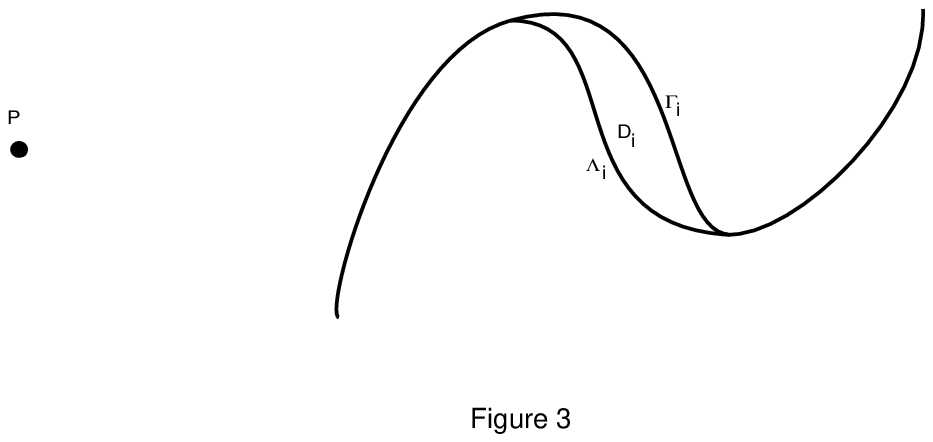}}$$
\endinsert
\par Now we wish to perform ``ambient singular surgery'' to $\cF$  along  $\Cal
M$ so as to obtain a floppy curve $\cf \subset \cp$. What we actually do is
perform a similar construction using a segment of the pencil of lines through
$P^*$ where $P^*$ is near $P$ but not in $\br \cP.$ Note that for each
inessential path $\L_i$,
we then obtain two paths $\G_i$ and $\G^{*}_i$ on $\rF$ and one of these will
be nested within the other. For notational convenience assume that  inner
segment is $\G^{*}_i.$ If  we choose our
inessential $\L^{*}_i$ so that $D_i^*$ is also nested within $D_i$, and we
choose a shorter length for the line
segment field over the $D_i^*,$ then our resulting membranes $M'$ and $M^{*'}$
will only intersect in $\rp$.
Moreover
$(M'\cup M^{*'})\cap\cF$ will be the boundary of a union of annuli. Thus we may
delete
the interior of these annuli
and adjoin $M'\cup M^{*'}$ to obtain an immersed surface $\eufm F$. This
surface may not be orientable.  However if it is orientable then it will be  a
floppy curve of degree m.
$\Bbb R \eufm F=\eufm F\cap \rp$ may be obtained by performing ``ambient
singular surgery''
to $\rF$ along $\cup_j \J_j$. In other words, we delete neighborhoods of the
end points of the essential arcs and attach $(\cup_j \J_j)
\cup (\cup_j \J^*_j)$. The line field $\Cal L_{\Bbb R \eufm F}$ along $(\cup_j
\J_j)$ is
given by $\Cal L_P$ and the line field
along  $(\cup_j \J^*_j)$ is given by $\Cal L_{P^*}$.
$\nu_\pm(\cf)=\nu_\pm(\cF)$.

\par Note that $\eufm F-\Bbb R \eufm F$ is oriented. Thus  $\eufm F$ will be
orientable if and only if  $\eufm o_\rF$ restricted to $\rF-$(neighborhoods of
essential arcs) extends to an orientation $\eufm o_{\eufm F}$. Below will given
an example where the orientation does not extend and the resulting $\eufm F$ is
nonorientable. One could study ``non-orientable'' floppy curves with the extra
data of a normal Euler number.
We do not pursue this now.

If $\eufm F$ is orientable we denote it by $\cf$.   $\cf$ will be dividing if
and only if $\cF$ is dividing and the complex orientation on
$\rF-$(neighborhoods of essential arcs) extends to an orientation on $\rf$.
\par
We may also perform a sequence of these singular surgeries to
$\cF$ along a collection of membranes which arise from a disjoint sequence of
pencil segments around the same point $P$.  There are probably situations of
interest where one can do a sequence of such modifications using membranes
arising from pencils around different points as well. But one would need to
know that
the portions of the membranes lying over disks near the inessential segments
do not intersect. For now, we will only consider  disjoint
segments of pencils around a single point.
\par We will present in Figures 4,5, and 6  three typical examples of how to
obtain a encoded picture of the $\rf$ from $\rF$.  We only draw relevant
portions of $\rF$
in a region around $\cP^+$. We draw in gray some of the lines in the pencil of
curves, as this helps us see the line fields $\Cal L_P$ and $\Cal L_P^*$. We
also draw the essential arcs with dashed lines in the ``before'' pictures.
A membrane like that used in Figure 4 will be called a simple membrane, and the
resulting surgery be called a simple surgery.
\midinsert
$$\vbox{\epsffile{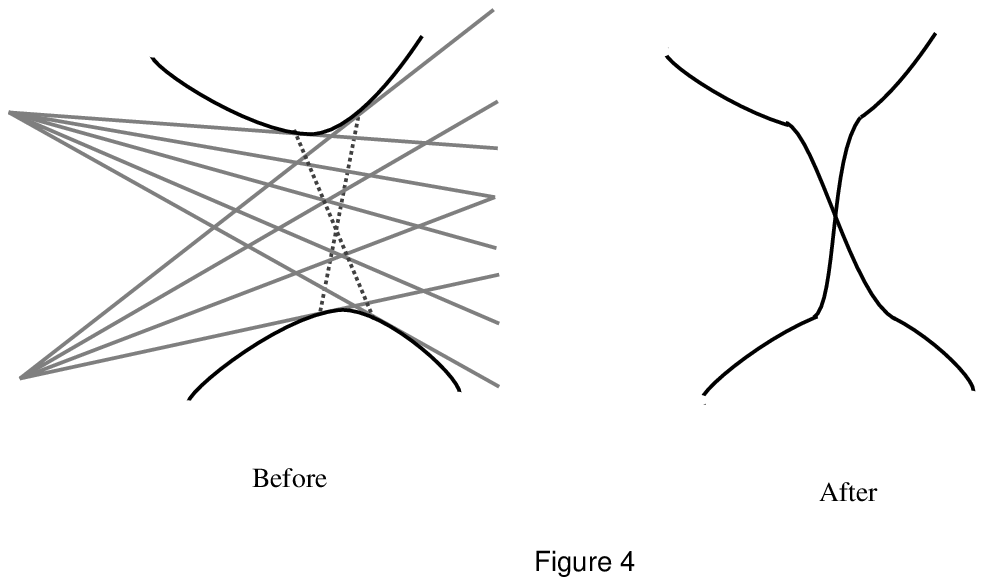}}$$
$$\vbox{\epsffile{fig05+.ai}}$$
\endinsert
Figure 5  has a thin part between the two flops. The ``figure eight'' in Figure
6 also has a thin part and a thickened part so the double point is a negative
double point.  Figures 4 and 6 use membranes of the type used by Viro and
Shustin. However the membrane used in Figure 5 is of a different type as a
essential arc crosses $\rf$. It is useful to know how these moves effect
$\D (\cf)$. A simple surgery leaves $\D (\cf)$ unchanged. The moves of Figures
5 and 6 lower $\D (\cf)$ by two. The move of Figure 6 will create a nondividing
curve from a dividing
curve. In Figure 7 we show an ambient singular surgery which leads to a
nonorientable curve.
\par

The example in Figure 4 could be a hyperbola that was constructed by perturbing
a pair of straight lines intersecting in  a single double point as in Gudkov.
Then our singular ambient surgery simply undoes the perturbation giving us a
pair of 2-spheres meeting transversely in point.
We may also topologically undo the move that replaces an isolated double point
with a small oval around that point
as follows.  Suppose that we have a floppy curve with an empty oval $C$, and
suppose
that $\Cal {L}_{\rf}$ is tangent to $C$.  Then we may  extend this to a  vector
field ${\eufm f}$ on the interior of $C$, and replace a tubular neighborhood of
$C$ in $\cf$ with the pair of disks traced out by the endpoints of $\pm i
{\eufm f}$. If the original floppy curve is dividing,  the new one  will be
dividing as well.  This is also a singular ambient surgery.
\midinsert
$$\vbox{\epsffile{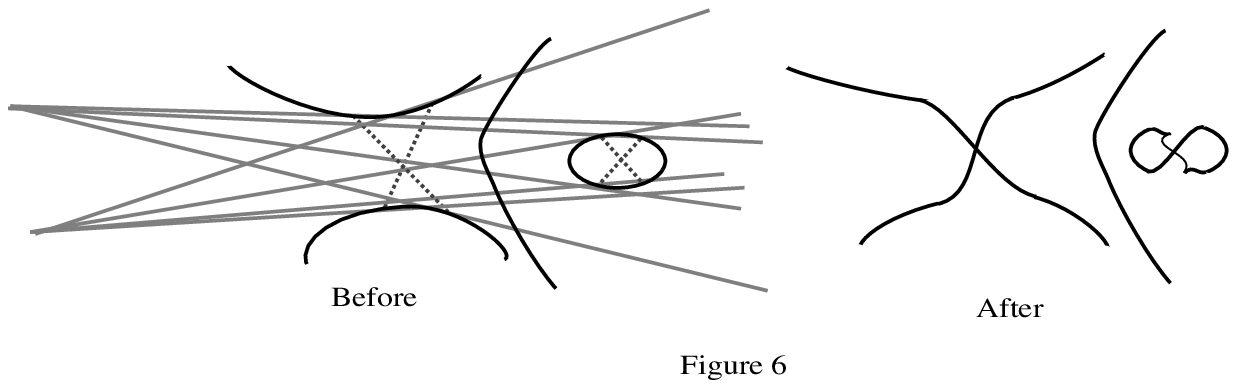}}$$
$$\vbox{\epsffile{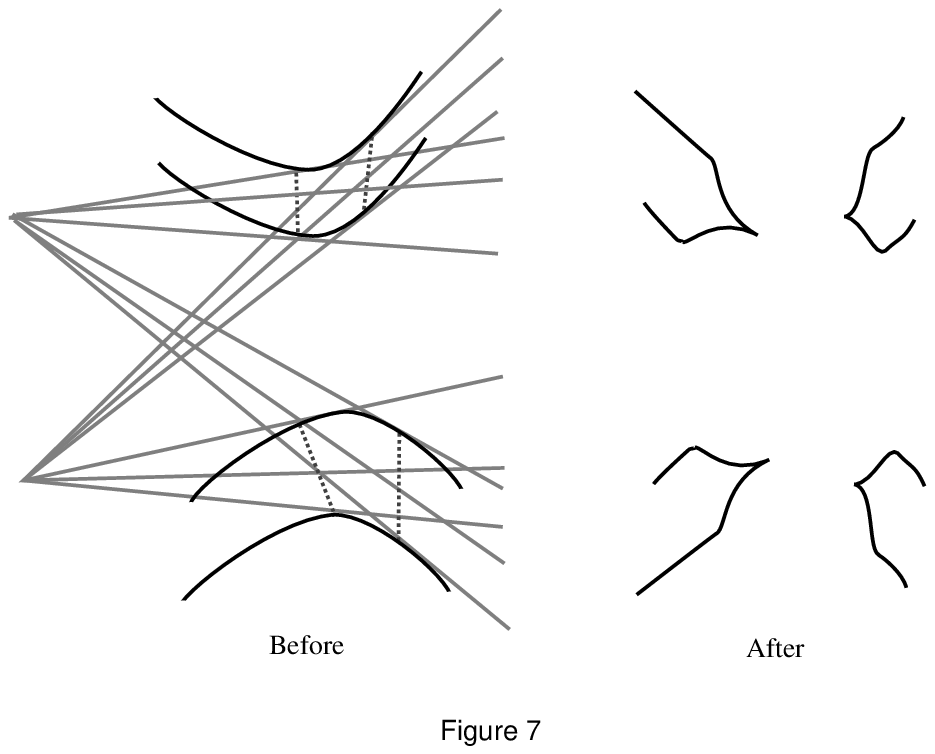}}$$
\endinsert

\head \S 3 Floppy curves and link cobordism\endhead

This section  generalizes \cite{G1} so that it applies to
floppy curves. The later sections do not depend on this section.  The reader
will find consulting \cite{G1} helpful.
Let $Q$ be the projective tangent bundle of $\rp$. A point in $Q$ is a pair
consisting of a point $p\in \rp$ and a line in $T_p(\rp)$. We let
 $\eufm f \in H_1(Q)$ be the class generated by the fiber. Since this is a
2-torsion class, we do not need to specify an orientation.
 \par
   Suppose that $C$ is the image of an immersion $\eufm j$ of a 1-manifold
$\hat C$
into $\rp$, and and $C$ comes equipped with a line field
$\Cal L_C$  which associates to $p \in \hat C$ a line in
$T_{\eufm j (p)}(\rp)$.  We assume that $\Cal L_C$ has the property that if
$p_1 \ne
p_2$, $\Cal L_C(p_1) \ne \Cal L_C(p_2)$. The image of a single circle is called
a constituent of $C$. Then we may define a link  $L(C)$ in $Q$ which lies over
$C$, it is the image of the map $\hat \eufm j$ from $\hat C$ to $Q$
which sends $p$ to the pair $(\eufm j (p), \Cal L_C p)$.  If $C$ has an
orientation we may lift the orientation to $L(C)$.
We may encode the isotopy type of the pair $(C,\Cal L_C)$ by a curve with flops
in the same way we encoded the isotopy type of the real part of a floppy curve.
 The isotopy type of $L(C)$ only depends on the isotopy type of $(C,\Cal L_C)$.
We will will also wish to consider $(C,\Cal L_C)$ enriched with an orientation
$\eufm o_{C}$ which assigns continuously for each $p \in C$  an orientation on
$j_*(T_p(C) \oplus (\Cal L_{C}(p))$ thought of as a plane in the tangent bundle
to the  the tangent disk bundle of $\rp$ at $(p,0)\in T(\rp)$. The isotopy type
of a triple $(C,\Cal L_C,\eufm o_{C})$ may be encoded by a curve with flops and
thickenings. The existence of an orientation $\eufm o$ implies that each
constituent must have an even number of flops.

\par
 We wish to specify a framing for $L(D)$.  For each orientation preserving
constituent  $\d$ of $D'$, we draw a parallel curve $\ddot \d$  whose only
intersections with $\d$ are due to double points of $\d$. For each orientation
reversing constituent $\d$ of $D'$, we draw a parallel curve $\ddot \d$ with
only one intersection with $\d$ not due to double points of $\d$. Let $\ddot D$
be the union of these pushoffs with line field $\Cal L_{\ddot D}$ obtained by
translating $\Cal L_{ D}$. The $L(\ddot D)$ is a push-off of $L(D)$ and so
gives a framing.
\par
We let $\tilde Q$ denote the tangent circle bundle of $\rp$.
 Let $\tilde L(C)$ be the inverse image of $L(C)$ under this two fold covering
map.   We can use $\eufm o$ to define an orientation on  $\tilde L(C)$ whether
or not $C$ has one. This is called the natural orientation. For this purpose it
is best to assume that our triple $(C,\Cal L_C,\eufm o_{C})$ has been isotoped
so it
is given by a curve with flops and thickenings.
Then for almost all points in $\tilde p \in\tilde L(C)$ are given by a point on
$p \in C$ together with a direction tangent to $C$ at $p$. This tangent
direction to $C$ at $p$ can be lifted to an orientation to $\tilde L(C)$ at
$\tilde p$.  If $C$ is thick at $p$ we choose this orientation at $\tilde p$.
Otherwise we choose the opposite orientation. If $C$ is oriented we may pick
out a sublink $L_+(C)$ of $\tilde L$ with one component covering each
constituent of $C$. One simply chooses the component of $\tilde L$ which
induces the given orientation on $C$ under projection. The isotopy classes of
all these links only depend on the isotopy of $C$ with its extra data (line
field and orientation $\eufm o$). We let $g$ denote the generator of
$H_1(\tilde Q)\approx Z_4$ represented by $L_+$ of a straight oriented line.
\par
  A small isotopy will make a floppy curve transverse to the Fermat quadric
$\Si$ in $\cp$ with equation $z_0^2+z_1^2+z_2^2 =0$. From now on we assume
$\cf$ is transverse to $\Si$.
Let $\eufm e(\cf)$ be half the number of negative intersections of $\cf$ with
$\Si$. So the number of intersections is $2 m+ 4 \eufm e(\cf)$. Let $w (\rf)$
denote the number of orientation reversing constituents of $\rf$.
Let $r(\rf)$ denote the number of real nonisolated  double points of $\rf^*$
and let $r_{\pm}(\rf)$ denote the number these points which are
positive/negative double points of $\rf^*$.  Let $\eufm L_k$ denote some
collection of $k$ thick lines in general position in $\rp$.  Let $\eufm
L_{m,\eufm e}$  denote some collection of $m+2 \eufm e$ lines with $m+\eufm e$
thick and $\eufm e$ thin. Let $L_k$, $L_{m,\eufm e}$,$\tilde L_k$, $\tilde
L_{m,\eufm e}$, $L_{+,k}$, $L_{+,m,\eufm e}$ denote respectively $L(\eufm
L_k)$, $L(\eufm L_{m,\eufm e})$,$\tilde L(\eufm L_k)$, $\tilde L(\eufm
L_{m,\eufm e})$, $L_+(\eufm L_k)$, $L_+(\eufm L_{m,\eufm e})$, where for the
last two links we give the lines arbitrary  orientations.  These links are
well-defined up to isotopy and do not depend on the specific choice of lines or
orientations.     The following theorem is derived by exactly the same method
as \cite{G1, Theorem (6.1)}. Let $\breve \rf$ denote the curve obtained from
$\rf$ by replacing every isolated point by a small oval. Recall that if $\rf$
is dividing, $\breve \rf$ has a complex orientation. A generalized cobordism is
a cobordism which allows ordinary double points \cite{G2}. It is the union of
generalized cobordisms which are the images connected surfaces with boundary.
These are called components of the cobordism.

\proclaim {Theorem (3.1)}   Given a floppy curve $\cf$ transverse to $\Si$
of degree $m$  with real part $\rf$   then  there is
an associated  generalized cobordism $G$ in $I\times Q$  from $L(\breve \rf)$
to $\Cal
L_{m, \eufm e(\cf)}$ with $\nu_+(\cf)$ positive double points and
$\nu_-(\cf)$ negative double points
such that
\roster
\item[1] $2 \chi (G)= \chi (\cf) + r(\cf) - \i(\cf)$
\item[2] One may orient $\tilde G$, the inverse image of
$G$ in $I \times \tilde Q,$ so that it becomes an oriented generalized
cobordism
from $\tilde L(\breve \rf)$ to $\tilde \Cal L_{m, \eufm e(\cf)}$ with their
natural orientations. This orientation is reversed by  the covering
transformation of  the cover $\tilde G \rightarrow G.$

\item[3] $G\circ G = (1/2)(m^2 -w(\rf))-r_+(\rf)+r_-(\rf)$.  Here $G$ is pushed
off itself in $I\times Q$ so that along
$L(A)$ and $L_m$, it has been pushed off with their given framings.
\item[4] If $[L_+(\breve \rf)]= m g$, then for any orientation of $L(\breve
\rf)$ and
some orientation of
$L_{m,\eufm e(\cf)}$, we have that $ \g(G)=0$ . If $[L_+(\breve \rf)]= (m+2)
g$, then for any orientation of $L(\breve \rf)$ and some orientation of
$L_{m,\eufm e(\cf)}$, $ \g(G)=\eufm f$.
\item[5] $\cf$ is dividing if and only
if $G$ is orientable.  If $G$ is orientable, an orientation of $G$
extends the orientation on $L(\breve \rf)$ induced by some complex orientation
of $\breve \rf$.
\item[6] There is a 1-1 correspondence between the $\br$-irreducible components
of $\cf$ and the components of $G$.
\endroster
Moreover if $G$ is
a generalized cobordism  in $I\times Q$  from $L(\breve \rf)$ to $\Cal
L_{m, \eufm e(\cf)}$  satisfying (1),(2),and (3),
then there is a floppy curve  $\cf$ transverse to $\Si$
of degree $m$  with real part $\rf$ with $G$ as it's associated  generalized
cobordism.
 \endproclaim
The invariant $\g(G)\in H_1(Q)$ of condition (4) is defined  in \cite{G2,\S 2}.
  Recall $\g(G)$One can calculate $[L_+(\breve \rf)]$ as follows. First one
resolves each double point of $\breve \rf$ as in  Figure 2a or 2b. Let C be the
resulting curve. Then $[L_+(\breve \rf)]= a g$ where $a$ is twice the number of
flops in $C$ plus twice the number of 2-sided components of $C$ plus the number
of thick one sided components minus the number of thin one-sided components.
Note there can be at most one one-sided component in $C$.

If one applied the  generalized Tristram-Murasugi inequalities  \cite {G2,8.3}
to the generalized cobordism $G$ given above
one would probably obtain the results which we describe in the next section. At
least this is what happens in the case $\cf$ is a nonsingular real algebraic
curve \cite{G3}. Similarly applying \cite {G1,8.5} to $\tilde G$ leads to a
generalization of the Viro-Zvonilov inequalities. We are working on deriving a
further generalization of the Viro-Zvonilov inequalities by generalizing the
method of sections six and later of this paper. In the case where $\nu =0$,
and $G$ is nearly planar (this means Proposition (1.2)
is nearly sharp), we may apply \cite {G2,Theorem (6.7)} on the Arf invariant of
links. For nonsingular curves, we obtain \cite{G3}, the Gudkov-Rokhlin
congruence and related
congruences.  Presumably
one should also be able to obtain some of the results of \cite{KV} in this way.
 In general one may apply the results of \cite{G2} to the link cobordisms
arising from components of the inverse image of $G$ in $I \times \check Q$
where $\check Q$ is any one of the five nontrivial covering spaces of $Q$.
\par
 \par In particular if $\cf$ is dividing, we have a generalized cobordism
$G_+$ in $I \times \tilde Q$ between
$L_+(\breve \rf)$ to $L_{+,m, \eufm e(\cf)}$. By the proof of (1.2), if
$h(\cf)=0,$ $n(\cf)=c(\cf),$ and $\nu(\cf)=0,$ then $G_+$ is planar.  We may
then apply \cite{G2,Theorem (6.8)} which asserts that the Arf invariant of
proper links (with extra data) is an invariant of planar cobordism. Recall that
$\tilde Q$ is oriented diffeomorphic to $L(4,3).$ Its first homology is
generated by a class $g$ which is represented by $L_1.$
If $L$ is a link in $\tilde Q,$ and $\g \in H_1(\tilde Q)$ such that $[L]=
2\g,$ and $q$ is a quadratic refinement of the linking form
$\ell$ on $H_1(\tilde Q),$ then $(L,\g,q)$ is proper if
for each component $K$ of $L,$ $q([K])+ \frac{lk(K,L-K)}{2} \equiv \ell([K],\g)
\pmod{1}.$ In \cite{G2,\S6}, we defined an  Arf invariant for proper triples
$(L,\g,q).$  It takes  values in $\bq/8\bz$. Let $q_{\frac{-1}{8}}$,
respectively $q_{\frac{3}{8}}$ denote the quadratic
refinement of the linking form on $H_1(\tilde Q)$ which sends  $g$ to ${-1/8}$,
respectively ${3/8}$.

\proclaim {Proposition (3.2)} The triples
$(L_+(\eufm L_{2k,\eufm e}),k g,q_{\frac{-1}{8}})$,$(L_+(\eufm L_{2k,\eufm
e}),(k+2) g,q_{\frac{3}{8}})$ are proper. The triples
$(L_+(\eufm L_{2k,\eufm e}),k g,q_{\frac{3}{8}})$,$(L_+(\eufm L_{2k,\eufm
e}),(k+2) g,q_{\frac{-1}{8}})$ are  not proper. Moreover we have
$$Arf( L_+(\eufm L_{2k,\eufm e}),k g,q_{\frac{-1}{8}})
=\frac{k^2}{2}=Arf( L_+(\eufm L_{2k,\eufm e}),(k+2)g,q_{\frac{3}{8}})+2$$
 \endproclaim

\demo{Proof} The results on properness follow easily from the definition.
By \cite{G2,(6.5)} we need only work out $Arf( L_+(\eufm L_{2k,\eufm e}),k
g,q_{\frac{-1}8}).$ There is a spanning surface $F_{2k,\eufm e}$ for
$L_{2k,\eufm e}$ consisting of $k+e$ annuli with $\g=kg.$ Each annuli is
given by a pencil of lines joining two distinct lines.  The Brown invariant
$\varphi^{F_{2k,\eufm e}}_{q_{\frac{-1}8}}$ must vanish \cite{G2,5.8}.
Thus we are left with $Arf( L_+(\eufm L_{2k,\eufm e}),k g,q_{\frac{-1}8})=
\frac 1 2 e(F_{2k,\eufm e})-\l(L_{2k,\eufm e}).$ One calculates that
$e(F_{2k,\eufm e})=\frac{k+e}{2},$ and $\l(L_{2k,\eufm e})=\frac{e-
k(2k-1)}{4}.$
\qed\enddemo
\proclaim {Theorem (3.3)} Let $\cf$ be a floppy curve of degree $2k.$  Suppose
that $h(\cf)=0,$ $n(\cf)=c(\cf),$ and $\nu(\cf)=0.$  If $\rf$ is equipped with
one of its
complex orientations and $(L_+(\rf),k g,q_{\frac{-1}{8}})$ is proper then
 $$Arf(L_+(\rf),k g,q_{\frac{-1}{8}})
=\frac{k^2}{2}= Arf(L_+(\rf),(k+2) g,q_{\frac{3}{8}})+2$$
\endproclaim
\demo{Proof} Note that by (1.2), $\cf$ is dividing. Also by the proof of
(1.2), $\tilde G$ is planar. Therefore the cobordism $G$ of (3.1) is planar as
well. By (6.8) of \cite{G2}, $Arf(L_+(\rf),k g,q_{\frac{-1}{8}})
= Arf( L_+(\eufm L_{2k,\eufm e}),k g,q_{\frac{-1}{8}}).$
\qed\enddemo
\par
When Theorem (3.3) is applied to nonsingular curves, one obtains Fiedler's
congruences \cite{F1}. This will be discussed in \cite {G4}.  One would expect
that Fiedler's congruences \cite{F2} for certain flexible nodal curves
could also be obtained but we have not verified this. Nevertheless the example
given below shows that
Theorem (3.3) applies in
some situations where \cite{F2} does not.
\proclaim {Theorem (3.4)(Shustin)} There is no nonsingular
real algebraic curve of degree eight with isotopy type
$1<2\coprod 1<18>>$.\endproclaim
\demo{Proof} We let $\ra$ denote such a curve and seek a contradiction. Pick a
point $P$ interior to one of the empty odd ovals. Consider the pencil of lines
through $P$ as it sweeps across the nonempty odd ovals.  As the lines in the
pencil pass from one empty even oval to the next either the other empty odd
oval is not encountered or it is.
In the first case one can find simple membranes joining consecutive ovals.  In
the second case one can find two  membranes of the type pictured in Figure 5
joining
the empty even ovals to the empty odd oval. Depending on where
the empty odd oval appears either  all the even empty ovals may be partitioned
into pairs of ovals joined by a simple membrane  (Case 1) or all but two of the
can paired in this way (Case 2).

In Case (1) by Fiedler's theorem \cite {F3} \cite{V2,(4.9)}, half of the empty
even  ovals are  oriented one way and half the other. By Rokhlin's formula
\cite{V2,3.13},  the only possible complex orientation for the remaining ovals
is as pictured in Figure 8a. This orientation may be excluded by Paris's
extremal property of Rokhlin's inequality \cite {Pa1,Pa2}. However to
illustrate
our methods we now prohibit this orientation using Theorem (3.3), just as we
did in a preliminary version of this paper written before Paris's work.
$$\vbox{\epsffile{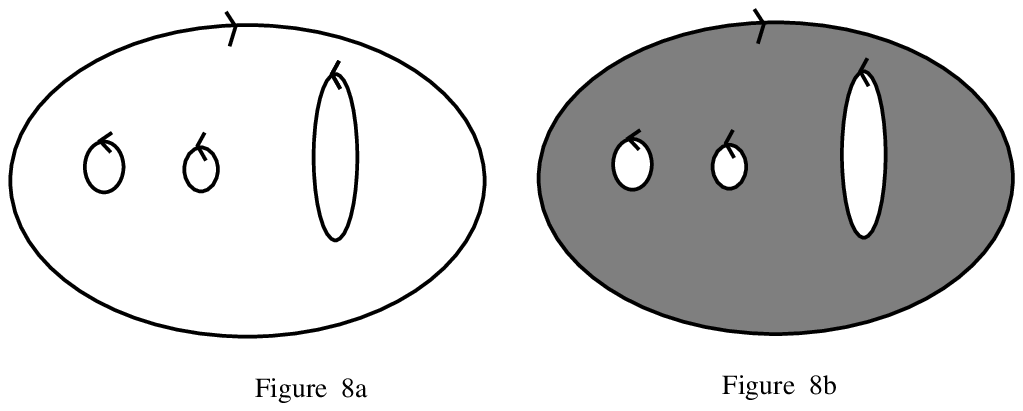}}$$
Perform singular ambient surgery along these membranes and obtain a floppy
curve $\cf$ with $\rf$ given by Figure 8a with nine additional figure eights
floating inside one of the odd ovals. Then $L_+(\rf)$ is given by this same
picture. But note that each component of $\rf$ which is given by a figure eight
is actually a small unknot unlinked from the rest of the link.
This is because when  traveling around a figure eight ones tangent vector does
not point in every direction. Thus the link given by a figure eight can be
described by the same diagram in a 3-ball and so bounds a disk disjoint from
the rest of the link. Let $L_1$ be the link in $\tilde Q$ which is defined by
Figure 8a. $L_1$ is proper for all $(\g ,q)$.  $L_+(\rf)$ and $L_1$ will have
the same Arf invariants for all $(\g, q)$.
 We can find  a vector field on the shaded region of Figure 8b tangent to
boundary along the boundary and pointing in the direction of  the orientation
of the boundary along the boundary which has  with a single zero of index minus
2 in the interior. Let $F_1$ denote this region blown up at this zero. In other
words a neighborhood of this zero is replaced by a Mobius band with core say
$\a$. As in \cite{G3}, \cite{G4}, the vector field then defines and embedding
of  $F$ into $\tilde Q$ with the boundary of $F$ given by
$L_1$. $\a$ is  a fiber of the bundle $\tilde Q \rightarrow \rp$ and so
represents $2g\in H_1(\tilde Q )$. We have that $\g(F) = 2g$, $e(F)= -4$,
$\l(L)= -3$.
$F$ capped of is just an $\rp$. One calculates that $\varphi
^{F_1}_{q_{\frac{3}{8}}} (\a)=  +1$.
See \cite{G2,\S6}, \cite{G4}. So $\b ( \varphi^{F_1}_{
 q_{\frac{3}{8} } } )=1.$  Thus $Arf(L_+(\rf),2
g,q_{\frac{3}{8}})=Arf(L_+(L_1),2 g,q_{\frac{3}{8}})= 2$.
This contradicts Theorem (3.3).

We now consider Case (2). We first observe that the arrangement in Figure 9a is
ruled out by Bezout's Theorem and a conic.
By Fiedler's theorem \cite {F3} \cite{V2,(4.9)}, the orientations of all but
two empty even ovals alternate and these two ovals must be oriented opposite to
the empty odd oval which does not include $P$. By Rokhlin's formula
\cite{V2,3.13},  there are only two only possible complex orientation for the
remaining ovals as pictured in Figures 9b and 9c. However the one shown in 9b
is shown to be impossible if we consider a pencil of curves  through a point
inside the other empty odd oval instead. Thus we have the complex orientation
as in Figure 9c.

$$\vbox{\epsffile{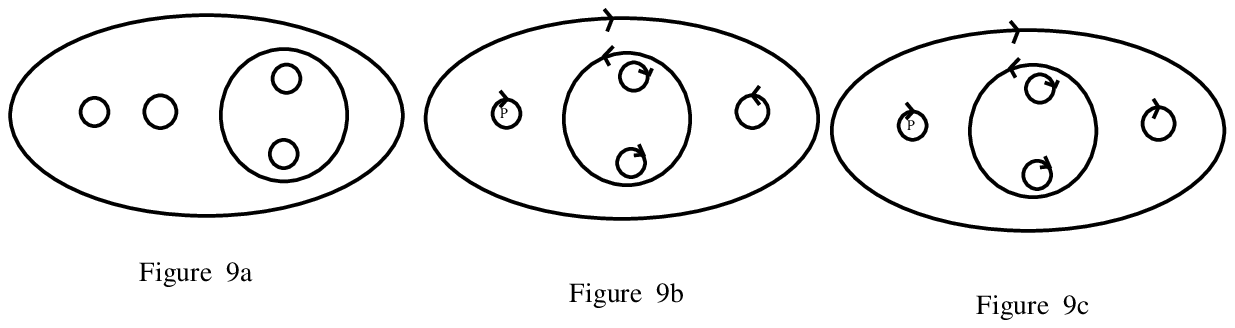}}$$
We perform eight simple surgeries and two of the type illustrated in Figure 5.
Let $\cf$ be the result. $\rf$ with eight figure eights deleted is pictured in
Figure 10a.  $L_+(\rf)$ is then the oriented link $L_2$ in $\tilde Q$ given by
Figure 10a together with eight unknots. 
$$\vbox{\epsffile{fig10+.ai}}$$
 Figure 10 illustrates a sequence of isotopic links in $\tilde Q$.  From
Figure 10d we see that $L_2$ is proper for all $(\g, q)$.  $L_+(\rf)$ and $L_2$
will have the same Arf invariants for all $(\g, q)$. Shade  Figure 10d
according to a `` checkerboard'' pattern, and then pick a vector field on the
shaded region tangent to boundary along the boundary and pointing in the
direction of  the orientation of the boundary along the boundary which has   a
single zero of index minus 2 in the interior.  As above this determines a
spanning surface $F_2$ for $L_2$. We have that $\g(F_2) =0,$  $e(F)=0,$
$\l(L_2)=-3,$ $\b (\varphi^{F}_{ q_{ \frac{-1}{8}} } )=1$. Thus
$Arf(L_+(\rf),0,q_{\frac{-1}{8}})= Arf(L_+(L_2),0 ,q_{\frac{-1}{8}})= 4$.
This contradicts Theorem (3.3).
\qed \enddemo

{\bf Remark} $L_+(1<2\coprod 1<18>> \text{(with any orientation)})$ is not
proper for any $(\g,q)$ because the components covering the empty even ovals
are improper.  A simple surgery and also the surgery of the type of Figure 5
change the  links in $Q$ and $\tilde Q$ corresponding to the real parts of the
floppy curves by band moves. The band moves associated to a simple surgery are
the band moves which undo the band moves  described in Figure 2a.  The band
moves associated to the surgery of the type of Figure 5 are described by the
line field with integral curves  pictured in Figure 11. In each case the
surgery to $\cf$ leads to a single band move to $L(\rf)$ and $L_+(\rf)$ and a
double band move to $\tilde L (\rf).$ Thus by \cite{G2, 6.10} the links $L_1$
and $L_2$ are guaranteed to be proper.  Thus ambient surgery sometimes allows
us to move from a situation where our results on Arf invariants do not apply
directly to one where they do apply.

$$\vbox{\epsffile{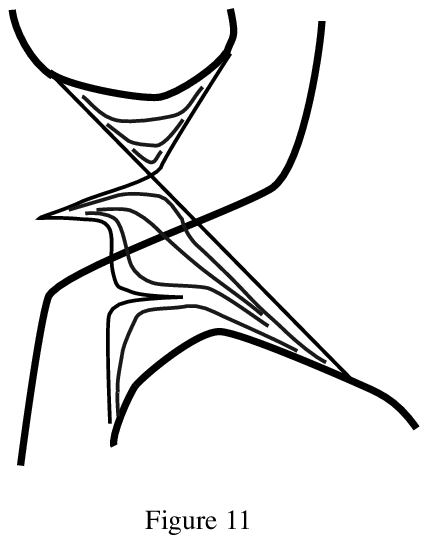}}$$

 \head \S 4 Generalized Viro-inequalities for floppy curves of even
degree \endhead

{\it For the rest of this paper we assume that the degree of $\cf$ is $2k,$ and
that $\cf$ is connected.}
In this section we generalize the inequalities that Viro originally found
for nodal curves so that they apply to floppy curves.  By a region  $R$ of
$\rf$ we
will mean a connected
region of $ {\Bbb R}P(2)- \rf^*$.   A point $c$ on the boundary
of $R$ which is a double point of $\rf$ is called a corner of
$R$.   A corner of $R$ is called simple if it has an open neighborhood
which intersects $R$ in a connected set.  If $ {\Bbb R}P(2)- {\Bbb
R}\cf$ is nonorientable, we let $w$ be a fixed orientation reversing
simple closed curve in $ {\Bbb R}P(2)- \rf$.
If $ {\Bbb R}P(2)- \rf$ is orientable, we pick  a fixed
orientation reversing simple closed curve $w$ on $\rf$.  A
non-simple corner of $R$ is said to cross $w$ if a neighborhood of
the corner intersected with the closure of $R$ is disconnected by its
intersection with $w$.

Let $s(R)$ be the number of simple corners of
$R$. If $k$ is even, we let $p(R)$ the number of non-simple  corners
of $R$ which do not cross $w$. If $k$ is odd we set $p(R)$ to be the
number of non-simple  corners of $R$. Let $\i(R)$ denote the number
of isolated singularities of $\rf$ in the interior of $R$. Let $O(R)$
denote the number of flops on the boundary of $R$ which point out of $R$ minus
the number of flops which point into $R$.  Let
$$\langle R, R\rangle = -2\chi (R) + {s(R) \over 2}+ 2 \ p(R)+ 2\
\i(R)+ O(R)$$
Note our regions $R$ are open. If $R$ and $R'$ are two distinct
regions of $B^{\pm}$ and  k is odd, let
$\langle R, R'\rangle$ denote one half the number of  points in the
intersection of the closures of $R$ and $R'$.  If k is even,
$\langle R, R'\rangle$ is one half the difference of the number of
intersection points which do not cross $w$  and  the number of
intersection points which do cross $w$. We call a region $R$
relevant if it is orientable or if $k$ is odd. Let $M^{\pm}$ denote the
matrix whose rows and columns are indexed by the relevant regions
of $B^{\pm}$ with $R,R'$
entry $\langle R, R'\rangle$.
\par
Given a matrix $M$, let $\sigma_ {\pm}(M)$ denote the number of positive,
respectively negative
eigenvalues of $M$. Let $\eta(M)$ denote the nullity of $M$. A vector
$\vec k$ with odd integer entries such that $M^{\pm}\cdot \vec k
=0$ is called a kernel vector for $M^{\pm}$.   Let $\e (\cf)$ be one if $\cf$
has a $\bc$-irreducible component of odd degree and be zero otherwise.

For a nodal curve $\rF$ one
lets $B^+$ and $B^-$ denote respectively the subsets of
${\Bbb R}P(2)$  where $F$ is greater than or equal to zero  and less
than or  equal to zero.   We  need to define analogous regions for a
floppy curve. Consider the new curve, say $\Cal R f$ consisting of simple
closed curves obtained if we ``resolve'' each double point of $\rf^*$
as in Figure 2.  We make an arbitrary choice at each crossing.
$\Cal R f$ will be null-homologous as well and thus it can contain
no 1-sided curve.   Now we can color  the components  of the
complement of  $\Cal R f$ in $\rp$ either black or white such
that each component of
$\Cal R f$  is on the boundary of exactly one black region.
There are exactly two  such colorings.  Then it is possible to pick a
related coloring of the complement of $\rf^*$ in a more or less
obvious way. Let $B^{*+}$ denote the closure of the black region and
$B^{*-}$ denote the closure of the white regions. Let $I^{\pm}$
denote the isolated points of $\rf$ in $B^{* \pm}$, and let
$B^{\pm}= B^{*\pm} \cup I^{\mp}$.  Let $O^{\pm}(\cf)$ denote the number of
flops in the line field which
point out of $B^{\pm}$ minus the number which point into $B^{\pm}$ . Note
$O^+(\cf)
= -O^-(\cf).$  Let $\i^{\pm}(\cf)$ denote the
cardinality of  $I^{\pm}(\cf)$.  Let $\i(\cf) = \i^+(\cf) +\i^-(\cf)$.
Of course if we
had made the opposite choice of original coloring, the roles of
$B^{\pm}$ would be reversed. Note that in (4.1$^{\pm}),$ (4.2$^{\pm}),$  and
(5.1$^{\pm}),$ the roles of these regions is completely symmetric.

{}From now on we will
abbreviate  $d_{\pm} (\cf)$ by $d_{\pm}$ , $\i^{\pm} (\cf)$ by
$\i^{\pm}$, $\D(\cf)$ by $\D$ etc.  Let
$$\Cal E^{\pm} = Max\{0,\eta(M^{\pm})-n + \e \}.$$

\proclaim {Theorem (4.1$^{\pm}$)} We have:
 $$ \sigma_ +(M^{\pm}) + \Cal E^{\pm} \le
(k -1)(k-2)/2 -\frac{\D-O^{\pm}}{4}\leqno(1^{\pm})$$
$$ \sigma_ -(M^{\pm}) + \Cal E^{\pm} \le
3k(k-1)/2 + \chi (B^{\pm}) - \frac{\D+O^{\pm}} {4} + \frac{r-\i-d_+ +d_-}{2}
\leqno(2^{\pm})$$
  \endproclaim
\par
For nonsingular curves Equations $(1^{\pm})$   are
the
strengthened Arnold inequalities and  Equations
$(2^{\pm})$  are the strengthened  Petrovski\u i inequalities. One
obtains from the  Addenda given below exactly the same extremal
properties for
these inequalities as are given in [W].

For nodal real algebraic curves, Theorem (4.1)$^{\pm}$ is almost the same as a
theorem  Viro announced in his abstract \cite{V3}, which he only states in
complete detail for $k\equiv 1 \pmod{2}$.  There is a misprint in the formula
which corresponds to Equation (2$^{\pm}$):  $-2\chi(\G)$ should read
$-2\chi(\bar \G).$
Our $\Cal E^{\pm}$ is sometimes a little bigger than the corresponding term in
\cite{V3}. Viro was aware that this could be improved. The left hand  side of
Equation (2$^{\pm}$) is smaller by $\i$ than the equation in \cite{V3}.
We improved our earlier version of Lemma (9.6) below after reading \cite{Fi}.  This
lead to an improvement in (4.1) so that our version of the Viro-inequalities
for nodal curves became the same as Finashin's.

We will say $\cf$ is 2-irreducible, if $\cf$ is either strongly-irreducible, or
if $n(\cf)=2$
and each $\bc$-irreducible component of $\cf$ has odd degree.  If $\cf$ is
2-irreducible these
inequalities can sometimes be improved.  The proof we give uses an adaptation
of a lemma of Viro and Zvonilov \cite{ZV} that we give
below as Lemma (7.3), as well as arguments from \cite{W}. Consider the
equations:
 $$ \sigma_ +(M^{\pm}) + \eta(M^{\pm}) \le
(k -1)(k-2)/2 -\frac{\D-O^{\pm}}{4}\leqno(3^{\pm})$$
$$ \sigma_ -(M^{\pm}) + \eta(M^{\pm}) \le
3k(k-1)/2 + \chi (B^{\pm}) - \frac{\D+O^{\pm}} {4} + \frac{r-\i-d_+ +d_-}{2}
\leqno(4^{\pm})$$

These are just Equations $(1^{\pm})$ and $(2^{\pm})$ with $\Cal E^{\pm}$
replaced
by $\eta(M^{\pm})$.
We have:
\proclaim {Addendum (4.2$^{\pm}$)} Suppose that $\cf$ is  2-irreducible,
then $(3^{\pm})$ and $(4^{\pm})$ can fail by at most one. If
either of
$(3^{\pm})$ or  $(4^{\pm})$  fail then all of the following conditions
hold:

(a) $\Bbb R \cf$ is dividing

(b$^{\pm}$) All of the regions of $B^{\pm}$ are relevant

(c$^{\pm}$) There is a kernel vector for $M^{\pm}$

\noindent Moreover, if $\cf$ is strongly-irreducible, and either of $(3^{\pm})$
or $(4^{\pm})$  fail to hold then
both of $(3^{\mp})$ and $(4^{\mp})$  will hold.
 \endproclaim

The hypothesis that $\cf$ is strongly-irreducible in the last line of (4.2) is
necessary as the example of two transverse real lines shows.

For the next addendum, we need to specify
precisely which region is $B^+,$ in the case where
every constituent on $\rf$ is null-homologous in
$\rp.$ Suppose  one assigns orientations arbitrarily to the constituents
 of $\rf.$  Choose $B^+$ so that $Int(B^+)$ consist of the set of points $x$
such that
$\rf^*$ is homologous to twice a generator for  $H_1(\rp - x).$ $B^+$ does not
depend on the arbitrary choice of orientations for the constituents.
Regions of the complement of $\rf^*$
which can be joined by  an arc which meets $\rf^*$ exactly once
transversely are assigned  opposite values by this procedure. Note that in the
case of a nonsingular curve $\rf$, this agrees with the usual convention than
$B^{-}$  should be non-orientable. Moreover if we smooth the intersections
compatibly with the orientation, then we will obtain a disjoint collection of
ovals. Then the $B^+$ for this collection of ovals determines a choice of $B^+$
for $\rf.$ This is the same $B^+$ defined above.

\proclaim {Addendum (4.3)} Suppose that $\cf$ is strongly-irreducible, every
immersed closed curve on $\ra$ is null-homologous in $\rp,$ and we
have chosen $B^+$ by the above scheme.
If either $3^+$ or $4^+$  does
not hold then $k$ is even.
 If either $3^-$ or $4^-$ does not hold then $k$ is odd.
\endproclaim

We also have the following proposition which follows from Equation (21) in \S
6, and the fact that each
constituent  of $\rf^*$ has an even  number of flops.

\proclaim {Proposition (4.4)} If $\cf$ is a floppy curve of degree
$2k$,  then $\D(\cf)\equiv O^{\pm} \pmod{4}$. Also $O^{\pm}$ is even.
So  $\D(\cf)$  is even, and  $\chi
(\cf)\equiv d(\cf) \pmod{2}.$ \endproclaim

Suppose we have found some simple membranes for a scheme for a hypothetical
complex curve $\Bbb CA.$ We can either follow
Viro's method of looking at the homology classes in the double branched cover
of $\Bbb CA$ given by the inverse image of the membrane, and then follow the
proofs of Arnold and  Petrovski\u i inequalities. Alternatively we could
perform simple surgery to obtain $\cf$ and apply (4.1) ,(4.2),(4.3) or (5.1) to
$\cf.$
Preliminary investigations indicate that one will get equivalent results either
way. Note this is only for simple membranes as in Figure 4.

\head \S 5 The Determinant of $2M^{\pm}$\endhead

The idea for the further conditions in this section came from the determinant
or discriminant conditions mentioned in \cite{S1,\S2} and \cite{KS,5.1}. Let
$\r^{\pm}(\cf)$ denote the number of relevant regions in $B^{\pm}$.
Let $\p^{\pm}(\cf)$ be the number of non-relevant regions in $B^{\pm}$. There
can be at most one such region.  Let
$$ h^{\pm}(\cf) = 3+r +
\chi(B^{\pm})+ \frac {d-\chi(\cf)}{2} -2\r^{\pm} -2\i^{\mp} - 2\p^{\pm}
,\leqno(5)$$
$$ P^{\pm}(\cf) =1+\chi (B^{\pm}) +\frac{r -\i-\chi(\cf)}{2}.\leqno(6)$$

\proclaim{Theorem(5.1$^{\pm}$)}   Let $\cf$ be  2-irreducible. We have that
$h^{\pm}(\cf)$ is a nonnegative integer.  If   $ \r^{\pm}(\cf)=P^{\pm}(\cf)$
then, $|\det(2 M^{\pm})|= b^2 2^{h^{\pm}+\r^{\pm} -d}$ for some integer $b$. \
\endproclaim

 The proof of this theorem  will be given in \S10.  The following observation
follows easily from the material in \S10:  in the situation where
$\eta(M^{\pm})=0$, one has that
$ \r^{\pm}(\cf)\le P^{\pm}(\cf)$ and moreover one has that: $ \r^{\pm}(\cf)=
P^{\pm}(\cf)$ if and only
if  Equations $(3^{\pm}$) and ($4^{\pm}$) are both sharp. The inequality
$h^-\ge 0$
when applied to a nonsingular real algebraic curve, where we choose $B^-$ to be
the nonorientable region, becomes Harnack's inequality. The only nonsingular
real algebraic curve where either of
$ \r^{\pm}(\cf)=P^{\pm}(\cf)$ holds is the degree two curve consisting of a
single oval.  The nodal real algebraic curve of degree  four constructed by
perturbing  some of the singularities  in a collection of four lines
illustrated in Figure 12 is a simple example of a curve where $
\r^{+}(\cf)=P^{+}(\cf)$ holds.  As yet we have no interesting prohibitions
arising from Theorem (5.1).
$$\vbox{\epsffile{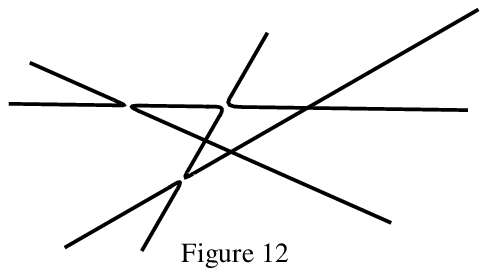}}$$

\head \S 6 The double branched cover of ${\Bbb C}P(2)$ along $\cf$
\endhead

Let $\cf_i$ be the $\bc$-irreducible components of $\cf,$ and $[\cf_i] = m_i
[{\Bbb C}P(1)]\in H_2({\Bbb C}P(2)).$ Let N denote a regular neighborhood of
$\cf$ in ${\Bbb C}P(2)$ which
is invariant under $T$.
By studying the long exact sequence of the pair, $({\Bbb C}P(2),{\Bbb
C}P(2)- N)$, one can show that there is a isomorphism from $H_1
({\Bbb C}P(2) - \cf )$ to $\Bbb Z^n/a\Bbb Z$ which sends the
meridians of $\cf_i$ to $e_i$ for each $i$, where $a = \sum
m_i e_i$.  Since, $\sum_i m_i$is even, there is a unique homomorphism $\phi
:H_1 ({\Bbb
C}P(2)
- \cf ) \rightarrow \Bbb Z_2$ sending the meridians of each
$\cf_i$ to the residue class of one. Note  $T_{|\cp - N}\circ \phi =
\phi$.
We now let $Y$ denote the associated branched cover of ${\Bbb
C}P(2)$.  $Y$ is a manifold except above the double points of $\cf$.
As the double branched cover of the Hopf link is
$RP(3)$\cite{R,p302},  the points above the double points have
neighborhoods which are cones on $RP(3)$.  So if $Y$ is a rational homology manifold and satisfies Poincare duality with rational coefficents.
 Let $\theta$ denote the
covering involution on $Y$ and $\pi$ denote the covering map from
$Y$ to ${\Bbb C}P(2)$.  There are two the involutions on $Y$  which
cover $T$. One will have $\pi^{-1}(B^{+})$ as fixedpoint set. The other
will have have $\pi^{-1}(B^{-})$ as fixedpoint set. We name the first
$T^+$ and the other $T^-$. Thus $T^-$ is the composition of $T^+$
and $\theta$.  Thus Y is a space with a $\Bbb Z_2 \times \Bbb Z_2$
action. Let $Y^{\pm}=\pi^{-1}(B^{\pm})$ , the fixedpoint set of
$T^{\pm}$. The orbit space of this action,
${\Bbb C}P(2)$ modulo complex conjugation is a smooth manifold
(actually the 4-sphere).   We can give a Whitney stratification of
this orbit space so that  the images of $\cf$ and ${\Bbb
R}P(2)$ are unions of strata. Thus we can triangulate the 4-sphere
so that the images of $\cf$ and ${\Bbb R}P(2)$ are sub-complexes
\cite{GM, p37 and references therein}.   If we then lift
this triangulation to $Y$, we have a simplicial $\Bbb Z_2 \times
\Bbb Z_2$ action on $Y$. By counting cells we see that

$$ \chi (Y^+) + \chi (Y^-) = 2 \chi({\Bbb R}P(2)) = 2. \leqno(7)$$

Because $\cf$ is
connected, the long exact sequence for the pair $(\cp,\cf)$ shows
that
$\b_1 (\cp,\cf,\bz_2)$ is zero. Thus by the sequence of Lee and Weintraub \cite{LW,Thm 1},$\b_1 (Y,\cf,\bz_2)$ is zero. 
Here we could as well use the Smith theory exact
sequence given in \cite{W, Appendix}.
The universal coefficient theorem then gives that $\beta _1(Y,\cf)$ is zero.  So by the long exact sequence of the pair $(Y,\cf),$  $\beta _1 (Y)$ is zero.
  By duality with rational coefficents  $\beta _3 (Y)$ also vanishes.

Let $_{\alpha} E_{\beta}$ denote the subspace of $H_2 (Y,\Bbb Q)$
consisting of those elements in both the $\alpha $ eigenspace of
$T^-$ and the $\beta$ eigenspace of $T^+ $, where $\alpha$, $\beta$
belong to the set $\{ -1,1\} $.  These subspaces are orthogonal with
respect to the intersection form on $H_2 (Y,\Bbb Q)$. Following
Wilson \cite{W}, we let $_{\alpha} E_{\beta}^\pm $ denote maximal
subspaces  of $_{\alpha } E_{\beta }$ on which the intersection form
is respectively positive  and negative definite.  The $+1$ eigenspace
for $\theta$ is $_{-1} E_{-1} \oplus _{+1} E_{+1}$ can be identified
by the transfer with $H_2({\Bbb C}P(2),\Bbb Q)$.  Since T acts by
multiplication by minus one on $H_2({\Bbb C}P(2),\Bbb Q)$ , the
intersection form is positive definite on $H_2({\Bbb C}P(2),\Bbb Q)$,
and $\beta_2({\Bbb C}P(2) = 1$, we conclude that $\dim( _{-1} E_{-
1}^+) =1$, and $\dim( _{-1} E_{-1}^-) = _{+1} E_{+1}^{\pm} =0$. Thus
$H_2 (Y,\Bbb Q)$ is written as a direct sum of five spaces:
$$  _{-1} E_{+1}^+ \ \oplus \  _{-1} E_{+1}^- \ \oplus \  _{+1} E_{-
1}^+  \ \oplus \  _{+1} E_{-1}^- \ \oplus \  _{-1} E_{-1}^+$$

Let $a$, $b$, $c$, and $d$ denote the dimensions of the first four of
these in order. Recall that the signature, $\text{Sign}(J)$ of an
involution $J$ on a rational homology 4-manifold $X$ is the
difference of the signature of the intersection form restricted to
the $+1$ eigenspace of $H_2(X,\Bbb Q)$ and the signature of the
intersection form restricted to the $-1$ eigenspace of $H_2(X,\Bbb
Q)$.  We also let $L(J)$ denote the Lefschetz number of $J$. We have:
$$ -(a + b) - (c + d) +1 +2= L(\theta)\leqno(8)$$
$$ -(a - b) -(c - d) +1 = \text{Sign}(\theta)  \leqno(9)$$
$$ (a + b) - (c + d) -1 +2= L(T^+)\leqno(10)$$
$$ (a - b) - (c - d) -1 = \text{Sign}(T^+).\leqno(11)$$

We now solve these equations obtaining:
$$a = (1/4)[ L(T^+) - L(\theta) + \text{Sign}(T^+)-
\text{Sign}(\theta)]+1\leqno(12)$$
$$b = (1/4)[ L(T^+) - L(\theta) - \text{Sign}(T^+)+
\text{Sign}(\theta)]\leqno(13)$$
$$c = (1/4)[ -L(T^+) - L(\theta) - \text{Sign}(T^+)-
\text{Sign}(\theta)]+1\leqno(14)$$
$$d = (1/4)[ -L(T^+) - L(\theta) + \text{Sign}(T^+)+
\text{Sign}(\theta)]+1.\leqno(15)$$

If $x$ and $y$ are 2-dimensional homology classes in a closed
rational homology 4-manifold, we denote their intersection number
by  $  x\circ y$. According to Theorem (7.1) and Lemma (7.2) in the next
section, we have that
$$\text{Sign}(\theta) = (1/2)[\cf]\circ [\cf] -d_+ +d_-  = m^2/2
-d_+ +d_-.\leqno(16)$$
By Lemma (9.5), we will show:
 $$\text{Sign}(T^+) = -\chi (Y^+)+ O^+.\leqno(17)$$
Using the oriented simplicial homology of Y \cite{Sp,p.159}, and the
Hopf trace formula \cite{Sp,p.195} , we see that:
 $$L(\theta) = \chi (\cf)\leqno(18)$$
 $$L(T^+) = \chi (Y^+).\leqno(19)$$
      By Equation (*) of \S 1 ,we have
 $$\chi (\cf) = 3m - m^2 + d_+ -d_-  + \D .\leqno(20)$$
 Substituting (16)--(19) in (11)--(15), and then substituting (20),(7)
and  $m = 2 k$ yields:
$$ a =  (k -1)(k-2)/2 - (\D - O^+) /4 \leqno(21)$$
$$ c =  (k -1)(k-2)/2 - (\D + O^+) /4 \leqno(21')$$
$$ b = 3k(k-1)/2 + (\chi (Y^+) - d_+ +d_-)/2- (\D+O^+) /4  \leqno(22)$$
$$ d =  3k(k-1)/2 + (\chi (Y^-) - d_+ +d_-)/2  - (\D-O^+) /4 \leqno(22')$$

Counting cells we see that:
$$ \chi (Y^{\pm}) = 2 \chi (B^{\pm}) +r(\rf)- \i(\rf) \leqno(23)$$
We will estimate the ranks of the $_{\alpha} E_{\beta}^\pm $ from
below using the homology classes residing on $Y^{\pm}$ to obtain
Theorems (4.1) and (4.2).

\head \S 7 Nice Involutions\endhead

In this section, we describe a class of  nice involutions, and
prove a
version of the G--Signature Theorem for them. We will see that the
involutions $\theta$ and $T^{\pm}$ on $Y$ are nice. Let $M'$ be a
compact smooth oriented $4$--manifold with boundary a disjoint
union of $RP(3)$'s. Let $M$ be the compact oriented rational
homology  $4$--manifold obtained by attaching the cone on $RP(3)$
to each boundary component.  Let $T'$ be an involution on $M'$ whose
fixed points set is a properly embedded (not necessarily orientable
or connected) surface $F'$ .  We assume that $T'$ acts on the boundary in
the following specified way. $T'$ may permute some pairs of
boundary components. Each of the boundary components invariant
under $T'$ is  oriented diffeomorphic as a manifold with action to
either the double branched cover of the Hopf link in $S^3$ , the
double covering of L(4,1) or the double covering of L(4,3). Then we
may extend $T'$ to an action $T$ on $M$  by coning off the action on
each boundary component which is not permuted and permuting the
cones on the boundary components which are permuted. We will call
any involution which can be constructed in this way nice.

   Each boundary component of $M'$ with a free involution will
contribute one isolated point to the fix point set of the involution on
 $M$.  Let $f_{\pm}$ denote respectively the number of boundary
components of $M'$ oriented as boundary of $M'$ with action the
double covering of L(4,1), respectively L(4,3).  Note the oriented
boundary of a neighborhood of a fixed point contributing to $f_{+}$ is
then the double covering of L(4,3). Let $F$ denote the fixedpoint set of
$T$ less the isolated fixed points set. We let $\bar M$ denote the
orbit space of the involution on $M$, and $\bar F$ denote the image
of $F$ in $\bar M$. $\bar M$ is also a rational homology manifold. Its
only singularities are cones on $RP(3)$ or $\pm L(4,1)$. $\bar F$ is
the  image of a immersed surfaces with only ordinary double point
singularities.  Note $\bar F$ lies completely  in the nonsingular part
of $\bar M$. We can view the pair ($M$,$T$) as the double branched
cover of $\bar M$ along the collection of immersed surfaces  $\{\bar
F_i\}_{i=1}^k$.  We let $e( \bar F )$ denote the Euler number of the
normal bundle of the immersion of $ \bar F$ in $\bar M$.

\proclaim {Theorem (7.1) }With $T$, $M$, and $\bar F_i$ as above:
$$\text{Sign}(T,M) =  f_+  -  f_-  + (1/2) e( \bar F ) .\leqno(24)$$
\endproclaim
\demo{Proof} If there are no singularities in $M$, this is just a
special case of the G--Signature Theorem \cite{AS}.  See \cite{JO}
for a topological proof of the G--Signature Theorem for involutions
or \cite{Go} for a topological proof of this theorem for 4-manifolds.

We deal first with the singularities of $M$ that are permuted in
pairs
by $T$. Let $P$ be some oriented smooth 4-manifold with boundary
$RP(3)$.  Replace neighborhoods  of each singular point in $M$ with
copies of P to obtain a new action .  It is easy to see the signature of
the new involution is the same as the old.   The left hand side of
(7.1)
is also unchanged. So it is enough to prove the theorem in the case
there are no singularities permuted in pairs.

Next we deal with singularities of $\bar M$ which are cones on
L(4,3).  We can replace a neighborhood of each such singularity with
 a disk bundle over a 2-sphere with Euler number four, and
simultaneously replace a neighborhood of the point above it with a
disk bundle over a 2-sphere with Euler number two.  The involution
extends over this new manifold.  The signature of the new involution
is one more than the old.  The right hand side of (7.1) changes
precisely the same way. We deal with singularities of $\bar M$
which
are cones on L(4,1) in a similar way: using a disk bundle with Euler
number minus four etc. Thus without loss of generality we may
assume that their are no isolated fixed points or singularities
permuted in pairs.

Form $M^-$ by deleting open neighborhoods of each singularities.  Let
$\bar M^-$ be the orbit space. The boundary of $\bar M^-$ is a
collection of 3-spheres. $\bar F$ meets each one in a Hopf link.
Form $\bar M^+$ by adding a pair of 2-handles to each 3-sphere along
these Hopf links, the first two handle one with framing one and the other with
framing minus one.  Then the branched cover of $\bar
M^-$ extends uniquely to a branched cover of $\bar M^+$ along $\bar
F^+$ which denotes  $\bar F \cap \bar M^-$ union the cores of the 2-
handles
with covering transformation $T^+$.
Above the branched cover of each pair of 2-handles,  we have a pair
of 2-handles which have been added to $M^-$ to form say $M^+.$  Each
time we add a pair of 2-handles we do not change the signature of
the intersection pairing restricted to the $-1$ eigenspace. Nor do we
change the signature of the pairing restricted to the $+1$
eigenspace or equivalently the signature of the orbit space by 1. This is
because the matrix $\bmatrix 1&1\\1&-1\endbmatrix$ has zero signature.
We also have $e( \bar F^+)=e( \bar F ).$
Thus we only need to show:
 $$\text{Sign}(T^+,M^+) = (1/2) e( \bar F^+ ) .\leqno(25)$$
The boundary of $\bar M^+$ consists of copies of $RP(3)$ and the
boundary of $ M^+$ consists a copy of the 3-sphere covering each of
the $RP(3)$'s. The cover of $RP(3)$ by $S^3$ has an orientation
reversing diffeomorphism, as $RP(3)$ does and the non-trivial cover
is unique.  Thus we may take two copies of $ M^+$ and glue them
together by an equivariant orientation reversing diffeomorphism on
the boundary.   In this way, we form $\tilde M$ with  involution
$\tilde T$.  By Novikov additivity,
$$\text{Sign}(\tilde T,\tilde M) = 2 \ \text{Sign}(T^+, M^+) .$$
$\tilde M$ is a smooth closed manifold and  $\tilde T$ is smooth.
We
know the theorem holds in this case. Thus twice Equation (25)
holds.\qed \enddemo

If we take the double branched cover of the 4-disk over the cone of
the Hopf link, we obtain a neighborhood above a double point in the
above construction. The branch set in this case is two 2-disks which
we can orient so that the double point is positive. Downstairs the
oriented  boundary of these two 2-disks is  the Hopf link with
linking number one.
Upstairs the oriented boundary of the lifts of theses two 2-disks is a
link in RP(3) with
linking number one half. To see this push of parallel copies of both
disks, they lift upstairs to  two 2-disks each lies in a neighborhood
of a disk of the branch set which double covers this disk.  These two
lifts  have  four positive transverse intersections. However the
boundary of each disk is two parallel copies of a component of the
boundary of the branch set.

This gives us a way to distinguish points lying above a positive
double point of $\bar G$ and points lying above negative double point of $\bar
G,$
when $\bar G$ is oriented.  If G intersects the RP(3) about the singularity in
a link with linking number $1/2$ ($-1/2)$, the double point below
the singularity is a positive (negative) double point of $\bar G$. In
this situation, we will refer the double point of G as a positive
(negative) singular double point.

The G-Signature Theorem for involutions is normally stated in terms
of the self intersection of the fixed point set.  We will need to apply
it in this form as well.   If we assume that $\bar G\subset \bar F$  is the
image of an
immersion of a closed orientable surface with ordinary double
points, then
we may distinguish between positive and negative double points. Let $G\subset
F$ be the corresponding subset of $F.$  We
let both $d( \bar  G) $ and  $d(   G ) $ denote the  number of double
points of $\bar G$ counted with sign.

\proclaim {Lemma (7.2)} For $\bar G$ orientable as above, we have:
$$   (1/2) e( \bar G) =   [G] \circ [ G]  - d(G ) = (1/2)[\bar G] \circ
[\bar G]  - d(\bar G ). \leqno(26)$$
\endproclaim
\demo{Proof}  $e( \bar G) $ is a signed count of the number of zeros
of a transverse section of the normal bundle of the immersion.  We
push off a parallel copy $\bar G'$ of $\bar G$ using such a section.
For each positive (negative) index zero of the section, we have one
positive (negative) point of intersection between $\bar G$ and $\bar
G'$.  Near each positive (negative) double point of the immersion  we
get an extra pair of positive (negative) intersections. Thus $[\bar G]
\circ  [\bar  G] = 2 d(\bar G) + e(\bar G)$.
Now if we consider  $ G'$ the lift of  $\bar G'$ in $M$.  We have  $
[G'] = 2 [G]$. On the other hand, we see that   $ [G'] \circ  [G]$ is also
$2 d(\bar G) + e(\bar G)$.  So
$$  [G ] \circ  [G] = (1/2) ( [\bar G] \circ  [\bar  G] ) \leqno(27). $$
The Lemma follows.\qed \enddemo

\par
If $\Cal M$ is surface (possibly non-orientable) with boundary (we
allow the boundary to have corners) mapped into $\bar M$ such that
$\partial \Cal M \subset \bar F$ , $\text{Int}\Cal M $ is immersed in $M$ minus
the singular set and $\text{Int}\Cal M$ is transverse to $F$, we call $\Cal M$
a membrane
for $\bar F$. Let $\tilde \Cal M$ denote the inverse image of $\Cal
M$ in $M$. Let $H_2(M,{\Bbb  Z}_2)^{T}$ denote the subgroup of
$H_2(M,{\Bbb  Z}_2)$ which is fixed by $T$. We have that $H_2(\tilde
\Cal M,{\Bbb  Z}_2)$  is one dimensional and is fixed by $T$, so
$\tilde \Cal M$ defines a class $[\tilde \Cal M] \in
H_2(M,{\Bbb  Z}_2)^{T}$. We will need the following lemma to prove
Addenda (4.2). Here we adapt an argument of Viro and Zvonilov's \cite{VZ \S2.2}
to our situation. We are working with a possibly nonorientable and possibly
singular branch set. Also $M$ may have singularities and is not necessarily a
$\Bbb Z_2$-homology manifold. However the first geometric argument in \cite{VZ}
does go through.
\par
\proclaim {Lemma (7.3) } There is a well-defined map $\Cal I
:H_3(M,{\Bbb  Z}_2) \rightarrow H_1(F,{\Bbb  Z}_2)$, given by
intersection.    There is a homomorphism $\Cal B$ from
$H_2(M,{\Bbb  Z}_2)^{T}$ to \break $H_1(F,{\Bbb  Z}_2)/ \text{Image
\ } \Cal I$. If $\Cal M$ is a membrane  for $\bar F$, $\Cal B[\tilde \Cal
M]=[\partial \Cal M]/{Image }\  \Cal I$.
\endproclaim
\demo{Proof}
 We will say a singular n-chain in $C_{n}(M,{\Bbb  Z}_2)$
is good if it is a sum of singular n-simplices of one of two kinds.  The first
kind has image in $M'$ , is smooth, and is transverse to $F$.  The second kind
is the join of a singular point with a singular (n-1)-simplex of the first kind
using the cone structure of neighborhoods.  We can pick a
good 3-cycle $\a$ representing a class $\bold a \in H_3(M,{\Bbb
Z}_2)$. The intersections of the images of its simplices with $F$
can be pieced together to form a 1-cycle in $F$. We define $\Cal I
(\bold a)$ to be the homology class represented by this 1-cycle.  If
$\a'$ is another good 3-cycle representing
$\bold a$, then $\a +\a'$
is a good cycle which must be the boundary of a good
4-chain $\b$.     Note that $\b$ intersected with $F$
provides a null-homology of $\a +\a'$ intersect $F$.
This shows that $\Cal I$ is well-defined.

If $\bold g \in H_2(M,{\Bbb  Z}_2)^{T}$, we can pick $\g$,  a good
cycle representing $\bold g$.
Then $\g+ T (\g)$ is the boundary of a good 3-chain $\d$. Then $\Cal
B (\bold g)$ is represented by the cycle obtained by piecing together
the intersection of the images of the simplices in $\d$ with $F$. If
$\d'$ is another 3-chain with boundary $\g+ T (\g)$, then $\d + \d'$ is
a 3-cycle whose intersection with $F$ defines an element in the
image of $\Cal I$. Thus the value of $\Cal B (\bold g)$ does not
depend on the choice of $\d$.\par

To see that the choice of $\g$ is not important, let $\g'$ denote a
second choice. Then $\g+\g'$ is the boundary of a good 3-chain $\z$.
Then $\d + \z + T (\z)$ is a 3-cycle with boundary $\g'+ T \g'$.  As
$\z$ and $T(\z)$ have the same intersection with $F$, $\Cal B$ is
well-defined. It is not hard to see $\Cal B[\hat \Cal M]=[\partial
\Cal M]/ \text{image } \Cal I$.\qed
\enddemo

\head \S 8 The  $\Bbb Z_2 \times \Bbb Z_2$ action above the real
double points of $\cf$\endhead

First we need to understand the involutions $T^{\pm}$ around the
singularities of $Y$ above the real double points of $\rf$.  For
this purpose we perform stereographic projection from a 3-sphere
$S$ in ${\Bbb C}P(2)$ surrounding such a point to see how the graph
$(\rf \cup {\Bbb R}P(2))\cap S)$ sits in $S$. This is an easy
exercise. We may take local coordinates so the point is the origin in $\Bbb
C^2.$ $\cf$ is locally given by the graph of $z_1 z_2 =0,$ and $T(z_1,z_2)=
(\bar z_1,\bar z_2).$ The result is shown in Figure 13a. The inverse image $R$
of
$S$ in $Y$ is thus the double branched cover of $S$ along a Hopf link
$H,$ and so is diffeomorphic to  $RP(3)$.  We first study the $\Bbb Z_2
\times \Bbb Z_2$ action on $R$.
\midinsert
$$\vbox{\epsffile{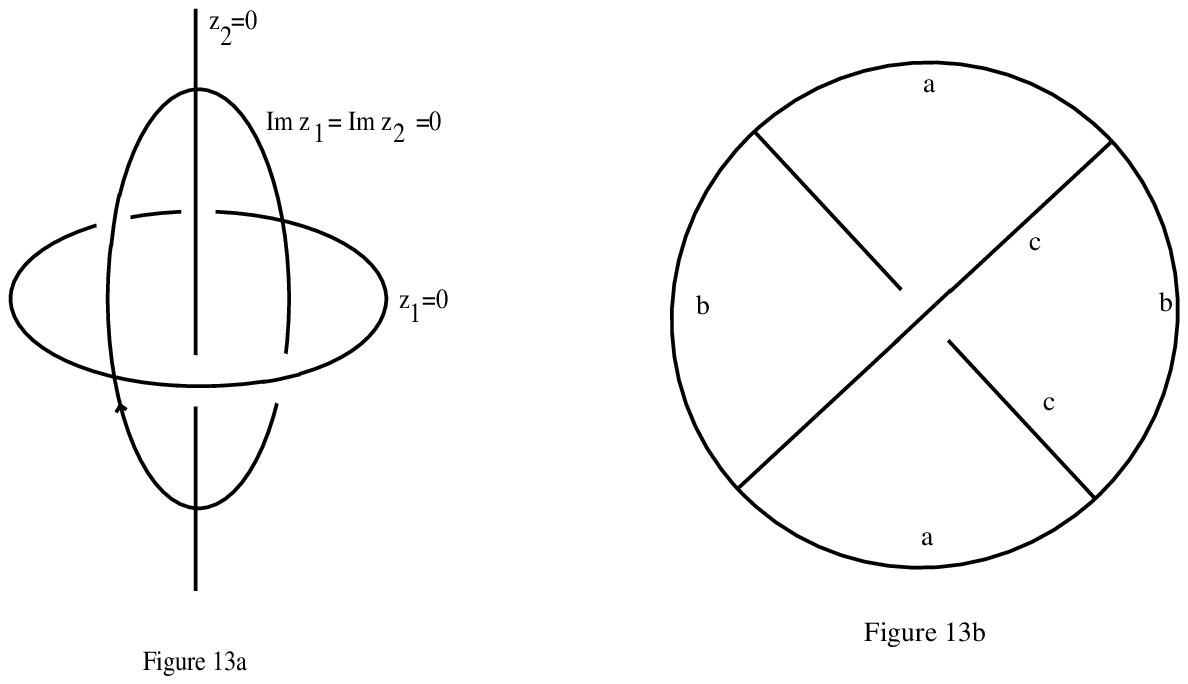}}$$
\endinsert
The orbit of the restriction of the involution $T$  to $S$  is a
3-sphere $\bar S$.
Figure 13b shows $\Gamma$ the image of $(\cf\cup {\Bbb R}P(2))\cap
S$  in  the orbit space of $T$.   The image
of  ${\Bbb R}P(2)\cap S$ consists of the edges labeled $a$ and $b$.
The image of  $\cf \cap S$ consists of the edges labeled $c$ .
The result is a completely symmetrical embedding of the complete
graph on four vertices.  In fact it is isotopic to the standard embedding of a
tetrahedron. There are both orientation preserving
diffeomorphisms and  orientation reversing  diffeomorphisms of the
3-sphere which send the graph to itself and performs any given
permutation of the vertices. To see this we make the following two
observations. An orientation reversing diffeomorphism which
performs any transposition is given by a reflection of the tetrahedron.  The
effect on the pictured graph of reflection  in the plane of the circle composed
of the edges labeled $a$ and $b$ can also be
achieved by an isotopy.

The inverse image of $S - ((\cf \cup {\Bbb R}P(2))\cap S)$ in
$Y$ is a regular $\Bbb Z_2 \times \Bbb Z_2$ cover of  the
complement of $\Gamma$ in the 3-sphere. This cover is classified
by a surjection  $\phi : H_1 (\bar S -\Gamma) \rightarrow \Bbb Z_2 \times
\Bbb Z_2$.  $R$ is the associated branched cover.  By
studying the relations coming from the 3--punctured spheres around
each vertex of $\Gamma$, we conclude that the $\phi$ must send the
meridian of an edge to its label where $a$, $b$,  are $c$  now the
distinct nonzero elements of $\Bbb Z_2 \times \Bbb Z_2$.  The
complete symmetry of embedding of $\Gamma$ implies that the
involutions $T^{\pm}$ and $\Theta$ restricted to  $R$ are all
conjugate.

We also need to know  certain linking numbers. $\rf \cap S$ consists of four
points. Let $x$ and $y$ be two
 points of $(\rf\cap S)$ adjacent on ${\Bbb R}P(2)\cap S$ and $\gamma$ be the
arc of
${\Bbb R}P(2)\cap S$ joining them. Consider a nonzero vector field
on a neighborhood of $\gamma$  in ${\Bbb R}P(2)$ which is tangent
to $\rf$ at $x$ and $y$ but points toward the central double
point at $x$ and away from it at $y$ and is only tangent to
$\gamma$ at exactly one point. If we multiply this vector field by $i$
it will be a normal vector field to ${\Bbb R}P(2)$ and therefore
tangent to $S$ along $\gamma$ . Let $\gamma '$ be an arc in $S$
with boundary
in $\cf$ obtained by pushing $\gamma$ off itself with this
vector field.
Let $\alpha$ and $\alpha '$ denote the inverse images of $\gamma$
and $\gamma '$ in the double branched cover of $S$ along $\rf
\cap S$. $\alpha$ is a simple closed curve in $R$ and $\alpha '$ is a
push off of $\alpha$.  Here the one parameter family of curves
upstairs traced out during the push-off upstairs covers the one
parameter family of arcs traced out by the push-off downstairs.  If
we assign an arbitrary orientation to $\alpha$ and the parallel
orientation of $\alpha '$, then the resulting linking number of
$\alpha$ and $\alpha '$ is well defined.
\proclaim{Lemma (8.1)} $Lk(\a,\a')= \frac{1}{2}.$
\endproclaim
\demo{Proof}  To see this we pick local coordinates so that
$\rf$ is the union of the coordinate axes, and take the
vector field to be tangent to the family of hyperbola given by $xy$ is
a constant. We parameterize $\gamma $  and $\gamma '$ and view
the images of the resulting curves under stereographic projection.
See Figure 14a.  In Figure 14b, we have  redrawn the  Hopf link given by
$z_1=0,$ and $z_2=0.$ We have also shaded the band between $\a$
and $\a'$ traced out by the one parameter family of arcs. We use the method in
\cite{R,p302} to see $\alpha$ and $\alpha '$ in $R.$  If we describe $R$ as
$-2$ framed surgery
to the unknot, then $\a$ is a meridian to the unknot and
$\a '$ is a parallel meridian.  It is not hard to see the linking
number is $1/2$.\qed\enddemo
\midinsert
$$\vbox{\epsffile{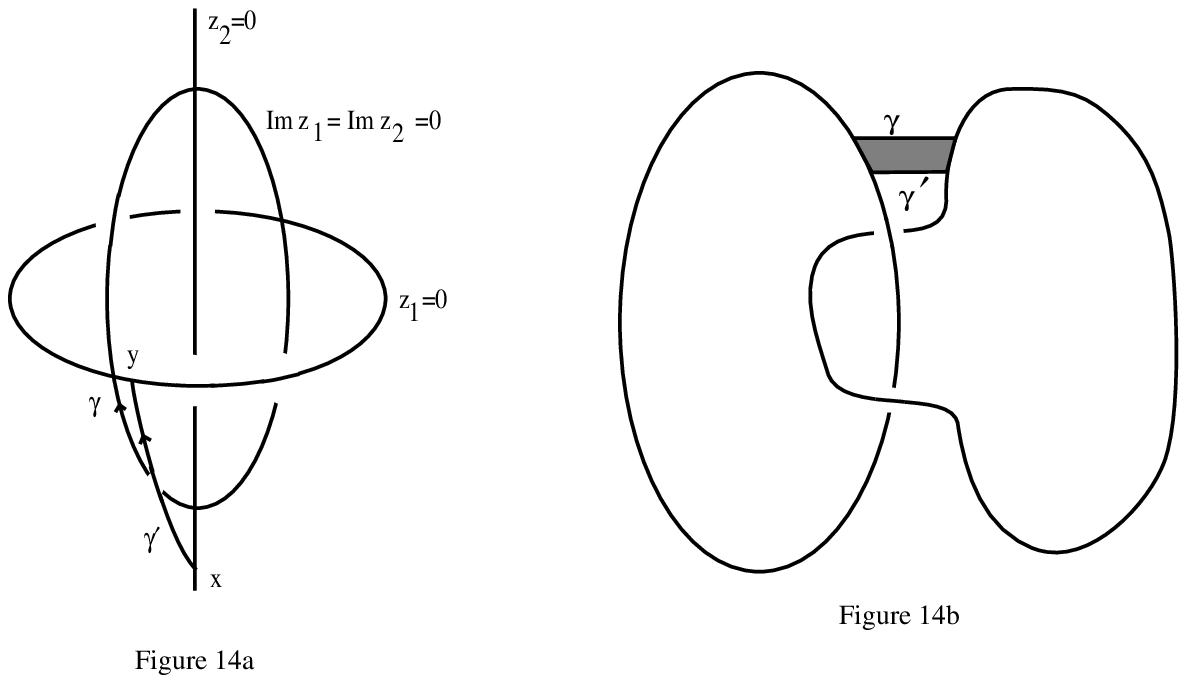}}$$
\endinsert
 \par
We will also need to know another linking number. Consider  a 4-ball $D$
neighborhood of our double point in $\cp$ with boundary $S$. In $D$
there is an annulus $A$ invariant under complex conjugation with boundary $\Cal
F
\cap S$ which one obtains by resolving the double point.
$A\cap \rp$ is a `` hyperbola''. Let $\hat N$ be the double branched cover of
$D$ along $A$. Consider the region $W \in N\cap \rp$ pictured in Figure 15.
\midinsert
$$\vbox{\epsffile{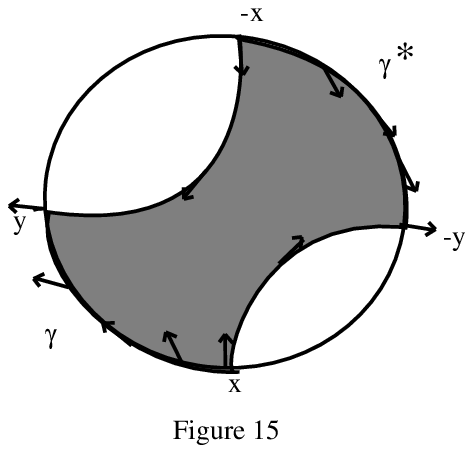}}$$
\endinsert
Let $\hat W$ be  the inverse image of $W$  in the branched cover. $\hat W$ is
an annulus. Give $\hat W$ an orientation then the boundary of $\hat W$ is an
oriented link with two components. One of these components is $\a$ lying over
$\g$. The other is say $\a^*$ lying over $\g^*$.

\proclaim{Lemma (8.2)} $Lk( \a,\a^*)= 1/2.$\endproclaim
\demo{Proof} Let $L$ be the link in $R$ consisting of $\a$ and $\a^*.$
 The curve labeled $\gamma '$ in Figure 14a can be
isotoped in $S -\gamma$ ($S$ is the boundary of $N$) keeping its
boundary on our branching set $S\cap \cf$ to the curve,
$\gamma^*$, given by reflection of $\gamma $ through the origin in
Figure 14a. Therefore $L$ is isotopic to the link consisting of  $\a$ and $\a
'$ with some orientation.  However we have
already seen that the linking number of the lift of $\a$ and
$\a '$ with (a push-off) orientation is $1/2$. So $Lk(\a,\a^*)= \pm
Lk(\a,\a^*)= \pm 1/2.$
\par
 We define a vector field $\vec  w$ on $W$ which is tangent to $A\cap \rp$
agrees with the vector field defined above along $\g$. Along  $\g^*$,  it is
obtained by reflection through the central point. We extend the vector field to
the interior of $W$ , so that it has a single zero of order
minus one.   We push $W$
off itself with $i\vec w$ to obtain $W'$. Let $\hat W'$
denote inverse image of $W'$ in the branched cover.  We orient
$\hat W'$  by pushing off the orientation on $\hat W$.  $L$ is the
boundary of $\hat W$, and we let $L'$ denote the boundary of $\hat W'$.  $\hat
W$
and $\hat W'$ will
intersect in two points above the zero of $\vec w$. Their
intersection number is $2$. See \cite{A,Lemma 6}, or
\cite{W, proof of (2.4)}.  The intersection form on $H_2(\hat N)$ is infinite
cyclic, and self
intersection of this generator is $-2$ \cite{V6},\cite{Ka}. So
$H_2(\hat N,\partial \hat N )$ is infinite cyclic and $[\hat W,
\partial \hat W]$ will represent some multiple of the generator say,
$\l$. It follows that
the total linking number between $L $ and $L'$ is
$2+ \l^2/2$. The formula we use here may be derived similarly to way one
derives the well-known formula which describes the $\Bbb Q/\Bbb Z$-linking form
of a 3-manifold given by surgery on a framed link with the inverse of the
linking matrix.  As this number is also $1+ 2 Lk(\a,\a^*),$ and $Lk(\a,\a^*) =
\pm1/2,$ we conclude $\l = 0$ and $Lk(\a,\a^*) =  1/2.$\qed \enddemo

 \midinsert
$$\vbox{\epsffile{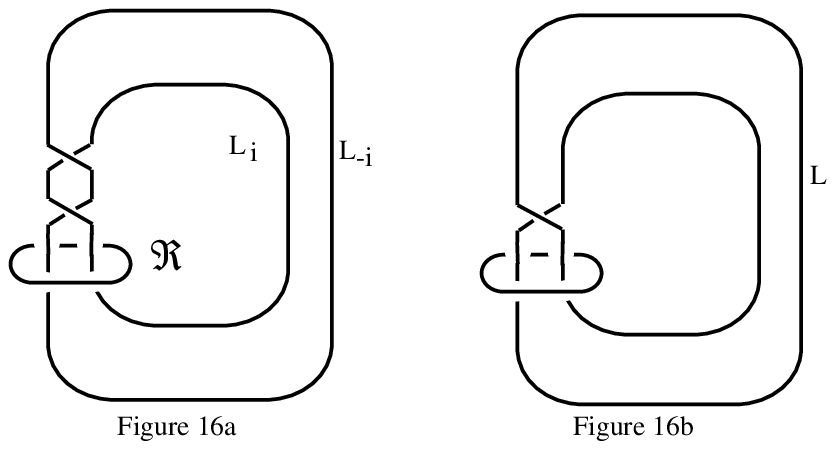}}$$
\endinsert

 We also must consider the involutions $T^{\pm}$ around the
singularities of $Y$ above the isolated  points of $\rf$.  In
\cite{G1,Figure 4a}, we have already described how $ \cf \cup
{\Bbb R}P(2)$ intersects  a small 3-sphere $S$ around such a point.
See Figure 16a. $\cf \cap S$
consists of the components labeled $L_{\pm i}$. $\eufm R$ is ${\Bbb R}P(2) \cap
S$.  We also described the image of
$(\cf\cup {\Bbb R}P(2)) \cap S$ in the orbit space $\bar S$ of
$T$ restricted to $S$ \cite{G1, Figure 4b}. See Figure 16b where we
have labeled the image of $L_{\pm i}$ simply $L$.  This link of two
components is symmetric in the sense that there is an isotopy which
switches the components. The link  $L_{+i} \cup L_{-i}$ is a Hopf link and
so the inverse image $\tilde S$ of $S$ in $Y$ is a copy of RP(3). We
have a $\Bbb Z_2 \times \Bbb Z_2$ action on $\tilde S$. The three
non-trivial elements $\Bbb Z_2 \times \Bbb Z_2$ each induce
involutions on $\tilde S$. The orbit spaces of these involutions must
be  the different 2-fold branched covers of $\bar S$ along one or
both the components of the link in Figure 16b. If we form the branched
cover along just one component we get $S^3$. Thus two of the orbit
spaces are 3-spheres.  The link in Figure 16b is isotopic to a $(2,-4)$
torus link and so using the method of Kirby and Akbulut \cite{AK}
the double cover is surgery on the unknot with framing -4 and so is
the lens space L(4,1).  Suppose that $p$  is an isolated double point in
$B^-$ and $\hat p = \pi ^{-1} (p)$. Then $T^+$ acts freely on the
boundary of a neighborhood of $\hat p $, and the orbit space of $T^+$
on this $RP(3)$ is $L(4,1)$. The orbit space of $T^-$ will be a 3-
sphere. If $p$
lies in $B^+$ instead the roles of $T^+$ and $T^-$ are
reversed.

Consider the component labeled $\eufm R$ in Figure 16a, we wish
to know the linking number of the two lifts of this component in
$\tilde S$. Viewing $L_{+i} \cup L_{-i}$ as a $(2,-2)$ torus link and using
the method of constructing double branched covers in \cite{R,p.302}
or \cite{AK}, we see that the lift of this component consists of two
oppositely oriented meridians in the description  of $\tilde S$ as
surgery on a unknot with framing $-2$. Thus their linking number is
$- 1/2$.

\head \S9  Proof of Theorem (4.1) and Addenda \endhead

Let $\phi :H_1 ({\Bbb C}P(2)- \cf ) \rightarrow \Bbb Z_2$ be as chosen in \S 6.

\proclaim{Lemma (9.1)} Let D be a oriented simple closed curve in $\rp-\rf$ .
If $D$ is one sided then $\phi( [D])=
\frac{m}{2}\pmod{2}$. If $D$ is 2-sided, $\phi( [D])=
0\pmod{2}.$ \endproclaim
\demo{Proof}   In the case $\cf$ is a nonsingular real algebraic curve, Lemma 1
of \cite{Pa1} provides an oriented 2-sphere $\Si \subset \Bbb CP(2)$ such that
$\Si \cap \Bbb RP(2)= D,$ $\Si \circ \cf= 2k,$ and $ \Si = D_1 \cup D_2,$
where $D_1$ and $D_2$ are disks, $D_1 \cap D_2 = D,$ and $T$  restricted to
$D_1$ is an orientation reversing homeomorphism to $D_2.$ The proof goes
through in this context. So $\phi( [D])= D_1 \circ \cf =k.$
\par
If $D$ is 2-sided, a similar proof gives a $\Si$ as above except $\Si \circ
\cf= 4k.$ Thus $\phi( [D])= D_1 \circ \cf =2k=0.$\qed\enddemo

   Recall that $Y^+$, the fix point set of $T^+$, is $\pi ^{-1}(B^+)$. If R
is a connected region of $B^+$, let $\hat R = \pi ^{-
1}\text{closure}(R)$.  Each $\hat R$ is the continuous image of a 2-
manifold
under a map which is one to one except at a finite set
where the map is two to one.  By an orientation on $\hat R$ , we
mean an orientation on its preimage.  If $R$ is orientable, then so is
$\hat R$ by Lemma(9.1), but an orientation on $R$ does not induce one on $\hat
R$.  We remark that even for nonsingular curves it not obvious that $\hat R$
will be orientable whenever $R$ is orientable. It is easy to construct an
involution on a Klein bottle with fix point set two disjoint circles with orbit
space an annulus.  See \cite {M,Remark (3.2)} for a discussion of this issue
for curves in $RP(1) \times RP(1)$. Also by Lemma (9.1), if $k$ is odd, $\hat
R$ is orientable whether or not $R$ is. For nonsingular curves these facts are
known \cite{W,(5.1)}.  To assign orientations to the orientable $\hat R$,
we consider the methods of removing double points discussed in \S 1.
Let $\cf_1$ denote the new floppy curve obtained in this manner if we
remove nonisolated double points by opening channels which join regions of
$B^+$,  and delete isolated double points. We can actually construct a
1-parameter family of floppy curves $\cf_t$ for $t \in [0,1]$ such that $\cf_t$
for $t \in (0,1]$ is an isotopy, and $\cf_0 =\cf$. We may then form $\check Y$
the double branched cover of $\cp \times [0,1]$ along $\cup_{t \in [0,1]}(\cf_t
\times t)$. Let
$S_t$ denote a region of $B_t ^+$,
which includes a region $R$ of $B^+$. Then we can view $\hat R$ as a
part of  the limit of $\hat S_t \times t$ as $t$ approaches zero.
\par
  We choose an orientation for each
orientable $\hat S_t$. If $\hat S_t$ is orientable and $R$ is included
in $S_t$, then the orientation on $\hat S_t$, induces one on $\hat R$.
Recall the fixed orientation reversing simple closed curve $w$
chosen in the introduction. If $\hat S_t$ is non-orientable then
$S_t$ will include $w$, and $\hat S_t - \pi ^{-1}(w)$ is orientable.
In this case choose an orientation on $\hat S_t - \pi ^{-1}(w)$ and
this will induce an orientation on $\hat R$ if $R$ is included in
$S_t-\{ w\}$. We have now chosen orientations on $\hat R$ for each
relevant region $R$ in $B^+$.

\proclaim {Lemma (9.2)}  Each $\hat R$ will have a positive singular
double point over each isolated point in $R$. $\hat R$ will also have
singular double points above the nonsimple corners of $R$. Such a
singular double point $c$ will be positive if k is odd or $c$ does not
cross $w$. Such a singular double point $c$ will be negative, if k is
even and $c$ does  cross $w$. Thus $\hat R$ will have $p(R) +\i(R)$
positive singular
double points. \endproclaim

\demo{Proof} The singular double points over each isolated points in
$R$ will all be positive by the last remark in \S 8.  To see this we
must note that an orientation on $\hat R$ will induce orientations on
the two sheets passing through $\hat p$ lying over an isolated point
$p$ which map down to opposite local orientations at $p$.

Note a neighborhood of a non-simple corner $p$, will intersect $R$ in a
region which looks locally like a neighborhood of the origin
intersected with the first and third quadrants in the plane. It is clear that
$\hat R$ will  have  a singular double point above such a corner of
$R$. It follows from our method of orienting $\hat R$ and Lemma (8.2)
that the sign of the double points will be as stated. \qed\enddemo

$Y^+$ is the union of all the $\hat R$ for $R$ in $B^+$ and $\i^-$
isolated points lying above $I^-$. $H_2(Y^+)$ is a free Abelian group
on the relevant regions of $B^+$.  We wish to figure out the
intersection form on this group.  Our choice of orientations, picks
out  a set of generators for this group : $[\hat R]$ as $R$ varies over
the relevant regions of $B^+$.

\proclaim {Lemma (9.3)} For R a relevant region of $B^+$, and any
orientation on $\hat R$, $ [\hat R]\circ [\hat R] = \langle R,
R\rangle$.\endproclaim

\demo{Proof} For simplicity we consider first the case when $O(R)$ is even. We
assume that we have already isotoped the floppy curve so that  so $\Cal
L_{\rf}$ is mostly tangent. We pick a vector field $ \vec v$ on R with isolated
zeros that is parallel to $\Cal L_{\rf}$. We also
model the vector field at the corners by the vector field $y{\vec i} -
x{\vec j}$ near the origin and in the first quadrant of the  $xy$
plane. We also insist that the zeros of $\vec v$ be simple in the
interior of $R$ and that there be a simple zero of order one at each of the
$\i(R)$ isolated double points of $\rf$ in $R$. The sum of the
$\pm 1$ indices of these zeros away from the isolated double points
is exactly $\chi (R)- \i(R)-O(R)/2$.  Then we can attempt to push $R$ off
itself in ${\Bbb C}P(2)$ with the vector field $i \vec v$. Let $R^+$
denote the result and $\hat R^+ = \pi ^{-1} R^+$.
We will succeed except at the isolated zeros of $\vec v$, which
include the corners.  $R^+$ and $R$ will represent the same homology
class in $Y^+$.  They will intersect each other transversely at
nonsingular points of $Y^+$ except above the corners  and  the
isolated singular points. The sum of the signs of these
intersection numbers will be $-\chi (R) + \i(R)+O(R)/2$. See \cite{A,Lemma
6}, or \cite{W,proof of (2.4)}. Thus the sum of the intersection
indices of $\hat R$ and $\hat R^+$ except above the corners and
isolated singularities  will be  $-2 \chi (R)+ 2\i(R)+O(R)$.  Above each
simple corner, the contribution of the intersection of $\hat R$ and $\hat
R^+$ to $[\hat R]\circ [\hat R]$ is
be computed as the linking number in $RP(3)$ of the curves
$\alpha$ and $\alpha '$ which we considered in Lemma (8.1).This is $1/2.$
\par
The
contribution above a non-simple corner $c$ is given by a linking
number in the copy of $RP(3)$ which is the boundary of a small
neighborhood of $\hat c$. $\hat R$ and $\hat R'$ will each intersect
this $RP(3)$ in a link of two components.
The contribution at a negative singular double point of $[\hat R]$ is
zero and the contribution at a positive singular double point of
$[\hat R]$ is two. See the proof of (8.2).
\par
 The contribution at an isolated singularity is zero.
This is calculated as follows. Push off with framing $-1$ a parallel
copy of the component labeled $\eufm R  $ in Figure 16a. Call it $\Cal
R^+$.  $\eufm R$ and $\eufm R^+$ describe the intersection of $R$ and
$R^+$ with a small 3-sphere $S$ centered at the singular point.
$L_{\pm}$ describes the intersection of $S$ with $\cf$.
$\hat S$, the branched cover along the link consisting of both
$L_{\pm}$, is the link of the singularity upstairs.  Our contribution
is the linking number of the inverse image of $\eufm R$ with inverse
image of $\eufm R^+$. However we must orient the two components of
each inverse image so that they induce opposite orientations
downstairs.  See the very end of \S 8. One checks that the linking
number is zero.
\par

For the case $O(R)$ odd, we need to replace $\vec v$ with  $\eufm l$, a line
field with isolated
singularities. We choose $\eufm l$ to be   parallel to $\Cal L_{\rf}$ and  at
the corners to be parallel to the vector field modeled on $y{\vec i} -
x{\vec j}$ near the origin and in the first quadrant of the  $xy$
plane. We also insist that there be a simple singularity of order two at each
of the
$\i(R)$ isolated double points of $\rf$ in $R$. Then we let $R^+$ denote the
double branched cover of $R$ along the singularities of
$\eufm l$ and corners of $R$, given by the endpoint of a small I-bundle
parallel to $i \eufm l$ completed by cones at the singularities of $\eufm l$.
Let $\tilde R^+$  be the inverse image of $R^+$ in $Y$. Then $\tilde R^+$
represents twice the class of $\tilde R$ and so the sum of the intersection
indices of $\tilde R$ and
$\tilde R^+$ is twice $ [\hat R]\circ [\hat R] $. The rest of the argument is
similar.
   \qed \enddemo

\proclaim {Lemma (9.4)} For $R$, and $R'$ distinct relevant regions
of $B^+$, with the given orientations on $\hat R$ and $R'$,$ [\hat
R]\circ [\hat R'] = \langle R, R'\rangle$.
\endproclaim
\demo{Proof} $\hat R$ and  $\hat R'$ intersects  only above
the points of intersection of closures of $R$, and $R'$. The
contribution there is given by the linking number we considered at
the end of the proof of Lemma (9.2).\qed \enddemo

Thus the matrix $M^+$ is the matrix which gives the intersection
form on $H_2(Y^+)$ induced from the intersection form on $H_2(Y)$
by the inclusion, with respect to the basis $[\hat R]$ as $R$ varies
over the relevant regions of $B^+$. Replacing $F$ by $-F$, and yields
the analogous result for $M^-.$

\proclaim {Lemma (9.5)$^{\pm}$} $\text{Sign}(T^{\pm})= -\chi(Y^{\pm})+
O^{\pm}$.
\endproclaim
\demo{Proof} By \S 8, we see that the action of $T^{\pm}$ on $Y$ is
nice in the sense of \S 7. Temporarily switching
the names of $B^{\pm}$ if necessary we may assume that every region of $B^+$ is
relevant. Recall that $Y^+$, the fix point set of $T^+$, is the union of
all the $\hat R$ for $R$ in $B^+$ and $\i^-$ isolated points. By
Theorem (7.1) and  the discussion in \S 8, each isolated point
contributes $-1$ to $\text{Sign}(T^+)$.  By Lemma (7.2), the
contribution of each $\hat R$ is $[\hat R] \circ [\hat R] - d(\hat R)$.
Let $n(R)$ be such that
$p(R) +n(R)$ is the number of non-simple corners of R.  Let $c(R)=
s(R) + 2p(R) + 2n(R)$. So $c(R)$ is the number of corners of $R$
where we count a corner twice if it is not
simple ( i.e if it occurs if twice while traveling around the
perimeter of $R$). By Lemma (9.2), $\hat R$ will have $p(R)+\i(R)$
positive singular
double points and $n(R)$ negative singular double points. So $d(\hat
R)= p(R)-n(R)+\i(R)$. On the other hand we can rewrite
$$[\hat R] \circ [\hat R] = \langle R, R'\rangle = -2\chi(R)+c(R)/2
+p(R)-n(R) + 2 \i(R)+ O(R).$$
Thus the contribution of $\hat R$ to $\text{Sign}(T^+)$ is $-
2\chi(R)+c(R)/2 + \i(R)+O(R)$. If we sum this contribution over all $R$ in
$B^+$, we will get
$-2\chi(Int( B^+))+r + \i^+ +O^+.$ So
$$\text{Sign}(T^+)= -2 \chi (\text{Int}( B^+)) + r + \i^+ -\i^- + O^+ = -\chi
(Y^+)+ O^+.$$

In the notation of \S 6, $\text{Sign}(T^-)$ is $(c-d) -(a-b)-1$. If we
plug in Equations (21),(21')(22),(22') and (7),
we see that $\text{Sign}(T^-)= -\chi(Y^-)+ O^-$. \qed \enddemo

The logic is slightly tricky here. We derive the equations of \S 6 for the
above choice of $B^+$ first, using Lemma (9.5)$^{+}.$ Then using the equations
of \S 6 for this choice, we derive Lemma (9.5)$^{-}.$ This may then be used to
derive the equations of \S 6 for the other choice of $B^+.$

\proclaim {Lemma (9.6)} The (-1)-eigenspace for the action of $\theta$ on
$H_3(Y,Y^+,\Bbb Q)$ has dimension less than or equal to $n-\e.$ \endproclaim

\demo{Proof}
Let $_{\alpha} {\eufm e}_{\beta}$ denote the dimension of the subspace of
$H_3 (Y,Y^+\Bbb Q)$
consisting of those elements in both the $\alpha $ eigenspace of
$T^-$ and the $\beta$ eigenspace of $T^+ $, where $\alpha$, $\beta$
belong to the set $\{ -1,1\}.$  Let $X_\pm$ denote the orbit space of the
involution $T^\pm$ on $Y$,
then as in \cite{Fi,4.3},  one may use the  Smith theory exact
sequence for  an
involution and the universal coefficient theorem to calculate that
$\b_3(X_+,Y^+) \le n-\e.$ Using the transfer \cite{B,2.4}, $(_+{\eufm e}_+) + (_-{\eufm e}_+) =
\b_3(X_+,Y^+).$
By the same argument, $(_+{\eufm e}_+) + (_+{\eufm e}_-) = \b_3(X_-,B^+).$ Using \cite{LW}, we have that
$\b_1(X_-)=0,$ as there is an involution on $X_-$ with $S^4$ as an orbit space with
a connected fixed point set.  By duality
$\b_3(X_-)=0,$ and so $\b_3(X_-,B^+)=0.$
Thus both $_+{\eufm e}_+$ and $_+{\eufm e}_-$ are zero.  So
$_-{\eufm e}_+ = \b_3(X_+,Y^+)\le n-\e.$ The eigenspace in question has dimension: 
$(_-{\eufm e}_+) + (_+{\eufm e}_-)\le   n -\e.$
\qed \enddemo

\proclaim {Proof of Theorem (4.1)}\endproclaim

Let $j^{+}$ denote the inclusion of $Y^{+}$ into $Y$. Then
$j^{+}_*$  maps  $H_2(Y^{+},\Bbb Q)$ into $_{-1} E_{+1}$. By Lemma
(9.6), the kernel of $j^{+}_*$ on $H_2(Y^{+},\Bbb Q)$ has dimension  less
than
$n -\e$.  It also has dimension less than $\eta (M^{+})$, as a
class which represents zero is in the radical of the intersection
form.
Pick a basis of eigenvectors for the matrix $M^+$. Let $S^+$ denote
the space spanned by the $\sigma_+(M^+)$ basis elements with
positive eigenvalue.  Let $S^0$ denote the space spanned by the
$\eta(M^+)$ basis elements with eigenvalue zero. Then the image of
$S^+ \oplus S^0$ in $H_2 (Y,\Bbb Q)$ will not intersect any of the
summands
$$\  _{-1} E_{+1}^- \ \oplus \  _{+1} E_{-1}^+  \ \oplus \  _{+1} E_{-
1}^- \ \oplus \  _{-1} E_{-1}^+.$$ So the inclusion of $S^+ \oplus S^0$
in $H_2 (Y,\Bbb Q)$ followed by the projection to $_{-1} E_{+1}^+$
will have kernel of dimension less than or equal to $\min \{ \eta
(M^+), n  -\e \}$. So we have:
$$ \sigma_+(M^+) + \eta (M^+)\le \dim( _{-1} E_{+1}^+)+ \min \{ \eta
(M^+),n -\e  \}.$$ Equation $(1^+)$ follows easily from this
and Equation (21).

Let $S^-$ denote
the space spanned by the $\sigma_+(M^+)$ basis elements with
positive eigenvalue. Considering the image of $S^- \oplus S^0$ and Equation
(22) as above
yields Equation $(2^+)$.
Equations $(1^-)$ and $(2^-)$ are obtained in exactly the same way, or we may
obtain them
by  reversing the roles of $B^+$ and $B^-$ by changing the sign of
the defining polynomial $F.$\qed \par

The following lemma follows easily from the Smith theory exact
sequence  or the exact sequence due to Lee and
Weintraub.

\proclaim {Lemma (9.7)}$\b_3 (Y,\Bbb Z_2) = n(\cf)-\e(\cf) -1.$  Thus  $\cf$
is  2-irreducible if and only if $\b_3 (Y,\Bbb Z_2) =0.$ \endproclaim

 \proclaim {Proof of Addendum (4.2)}\endproclaim

By (9.7), $H_3(Y,\bz_2) =0.$
By Lemma (7.3), the only possible nonzero
element of the kernel of  $j^{{\pm}}_*:H_2(Y^+,\Bbb Z_2)\rightarrow
H_2(Y,\Bbb Z_2)$ is
$[Y^{\pm}]$, and this  element can not be in the kernel
unless condition (a) holds. Here
$[Y^{\pm}]\in H_2(Y^{\pm},\bz_2)$ is the class of a cycle given the
sum of all the simplices in some triangulation of $Y^{\pm}$.

Thus $\b_3(Y,Y^\pm,\bz_2) \le 1.$
By the proof of (4.1$^\pm$), Equations (3$^{\pm}$) and (4$^{\pm}$) can fail by
at most $\b_3(Y,Y^\pm,\bq).$
By the universal coefficient theorem, $\b_3(Y,Y^\pm,\bq)\le \b_3(Y,Y^\pm,\bz_2)
\le 1.$ Thus Equations (3$^{\pm}$) and (4$^{\pm}$) can fail by at most one, and
then (a) holds.
\par

We now show that
the failure of  condition (b$^{\pm}$)
implies that $j^{\pm}_*:H_2(Y^{\pm},\Bbb Q)\rightarrow H_2(Y,\Bbb Q)$  is
injective. Let $\check Y$ denote the inverse image in $Y$ of the
union of the
relevant regions in $B^{\pm}$.  Then the inclusion induces an
isomorphism from
$H_2(\check Y,\Bbb Q)$ to $H_2(Y^{\pm},\Bbb Q)$. Thus it suffices to
show that
$H_3(Y,\check Y,\Bbb Q)$ is zero. This will hold if $H_3(Y,\check
Y,\Bbb Z_2)$
is zero. The boundary of the closure the union of any collection of
relevant
regions in $B+$ can not be null-homologous in $\cf$. So by Lemma
(7.3),
$j^{+}_*:H_2(\check Y,\Bbb Z_2)\rightarrow H_2(Y,\Bbb Z_2)$  is
injective. It
follows that $H_3(Y,\check Y,\Bbb Z_2)$ is zero.\par

We now assume that (3$^{\pm}$) or (4$^{\pm}$) fails, and wish to
establish
condition (c$^{\pm}$), having already established conditions (a) and
(b$^{\pm}$). $H_2(Y^{\pm})$
is free Abelian generated by $[\hat R_i],$ where $R_i$ ranges over all
the
regions of $B^+,$ each of which must be relevant.   It
follows that $j^{{\pm}}_*:H_2(Y^+,\Bbb Z)\rightarrow
H_2(Y,\Bbb Z)/  \text{torsion}$ is not injective.  So there must be an
element in  the kernel of the form $\sum_i k_i [\hat R_i]$
where the integral vector $\vec k=\{k_i\}$ is primitive. In particular some
$k_i$ is odd. Also we have that $M^{\pm}\vec k=0.$

$H_2(Y^{\pm},\bz_2)$ also has as a basis $\{[\hat R_i]\}$
(reduced  modulo two).
As $H_3(Y,\Bbb Z_2)$ is zero, $H^3(Y,\Bbb
Z_2)=0$,  and so $H_2(Y,\Bbb Z)$ has no 2-torsion. Thus the map
given by  reduction modulo two from $H_2(Y,\Bbb Z)$ to $H_2(Y,\Bbb
Z_2)$  factors through $H_2(Y,\Bbb Z)/ \text{Torsion}$.
It follows that  $\sum_i k_i [\hat R_i]$ reduced modulo two is in the kernel of
$j^{{\pm}}_*:H_2(Y^+,\Bbb Z_2)\rightarrow
H_2(Y,\Bbb Z_2).$ It is nonzero, since one $k_i$ is odd. In the proof of (a),
we saw
the only possible such element is $[Y^{\pm}]=\sum_i k_i [\hat R_i].$
It follows that all $k_i$ are odd, and $\vec k$
is a  kernel vector for $M^{\pm}$.  So
condition (c$^+$) holds.\par
To see the second statement, we note that $\rf$ has atleast one constituent as $\cf$ is strongly irreducible and must be dividing.
Then an argument in the
proof of [W,7.4] shows that not both of
$[Y^+]$ and $[Y^-]$ can be zero
in $H_2(Y,\bz_2)$. \qed
\par
\proclaim {Proof of Addendum (4.3)}\endproclaim

With our choice of $B^+$ the proofs of Lemmas (6.5), (6.6) and (6.7)
of [W] all go through. Then the proof [W, 7.4 (iii)] completes the
proof. \qed

\head \S10  Derivation of the Determinant Condition   \endhead

\proclaim {Proof of Theorem (5.1$^{\pm}$)}\endproclaim
We only need to prove the plus version. We desire to work with an actual
manifold rather than a rational homology
manifold so that we have a unimodular intersection pairing over $\bz$. So we
equivariantly resolve
$d$ singularities in the space with involution $(Y,T^+)$. Each singularity is a
cone on a copy of $\rpt$. Let $\Bbb D$ denote the disk bundle over $S^2$ with
Euler
class $-2$.

For each of the $\nu$  pairs of complex conjugate double points, we have a pair
of cones on $\rpt$ which are interchanged by $T^+$. We replace each of these
pairs with two copies of $\Bbb D$ which are interchanged by $T^+$. The sum of
their
fundamental classes is fixed by $T^+$.

Above each of the $i^- $ isolated double points in $B^-$,
we have an isolated fixed point of $T^+$ whose neighborhood is a cone on a copy
of $\rpt$ with a free action. We can equivariantly replace this neighborhood
with a copy of $\Bbb D$ with the involution which is the antipodal map on the
disk
fibers. The core 2-sphere is fixed by the involution, so its fundamental class
is also fixed.

Above each of the $i^+ $ isolated double points in $B^+$ and above each of the
$r$ nonisolated real double points,
we have a fixed point of $T^+$ whose neighborhood with involution is a copy of
the double branched cover of $B^4$ along
the cone on the Hopf link. We can replace this with the double branched cover
of  the 2-disk bundle over the 2-sphere with Euler class $-1$ along two
disjoint disk fibers. This is just another copy of $\Bbb D$
with a different involution. This changes the fixed point set by locally
replacing
the cone on two circles with the disjoint union of two disks. The homology
class of the core 2-sphere is fixed by the involution.

Let  $(\cy,\ct)$ denote the manifold with involution we obtain if we resolve
the singularities of $(Y, T^+)$ as above. As $\cf$ is 2-irreducible,
$H_3(Y,\bz_2)=0$. Comparing the Mayer-Vietoris sequence
for $Y$ as the union on the complement of  singularity set and a neighborhood
of the singularity set and a similar Mayer-Vietoris sequence for $\cy$, shows
that $H_3(\cy,\bz_2)=0$ and thus, using Poincare duality in an actual manifold,
$H_1(\cy,\bz_2)=0$. Thus $H_2(\cy)$ can only have odd torsion.
We have $$\align \b_2(\cy) &= \b_2(Y)+d= \chi(Y)-2 +d\\
 &=2 \chi( \Bbb CP(2))-\chi(\cf)-2+d= 4-\chi(\cf)+d.\endalign$$
 \par
The fix point set $\ef$ of $\ct$ can be obtained from $Y^+$ by replacing  $i^-$
isolated points with 2-spheres and by replacing $r+ i^+$ cones on circles with
pairs of disks. Thus we can find a basis $\Bbb B$ for
$H_2(\ef)$ indexed by the relevant regions of $B^+$ and the isolated points in
$B^-$. Thus $\b_2(\ef)= \rho^+ +i^-$ and $\b_2(\ef,\bz_2)=\b_0(\ef)= \rho^+
+i^-+ \p^+$.  We have that
$\chi(\ef) = \chi(Y^+) + i^- +r+i^+$. Then by Equation (23) we have:
$$\b_1(\ef,\bz_2)= 2 \b_0(\ef)-\chi(\ef)= 2(\r^+ +\p^+ +i^- -r- \chi(B^+)).$$
\par By a theorem of Floyd \cite{Fl,4.4}, \cite{B,p.126},
\cite{W,A2},$\sum_k \b_k(\ef,\bz_2) \le \sum_k \b_k(\cy,\bz_2)$. Equivalently
their difference:  $2- 2\b_0(\ef)+\b_2(\cy)-\b_1(\ef,\bz_2)$ is positive. By
the above equations, this last quantity is $2 h^+$. Thus $h^+$ is nonnegative.
The integrality of $h^+$  follows from Proposition (4.4).
\par
Let $H_+$ denote the
kernel of the endomorphism of $H_2(\cy)/\text{Torsion}$ defined by $1- \ct_*$.
Let $L$ be the subspace of $H_+$ generated by the images of elements of $\Bbb
B$   in $H_2(\cy)$, together with the $r+i^+$ classes coming from the 2-spheres
with involution replacing the singularities above the nonisolated real
points and the isolated points in $B^+$  and the $\nu/2$ classes given by the
sums of
paired of 2-spheres arising from the pairs of complex conjugate double points.
 $\dim (L) \le \Bbb L$ where we let $\Bbb L$ denote $\rho^+ + r+ i+ \nu/2.$ By
\S6,  subspace of on $ H_2(Y,\bq)$ fixed by the action of $T^+$ on $
H_2(Y,\bq)$ has dimension $a+b= P^+$. By
our hypothesis, this is also $\rho^+$. So $\dim (H)$ equals $\Bbb L.$ If
$\eta(M^+)$ is nonzero, $\det( M^{\pm})=0$ and the conclusion holds. Thus we
may assume that $\eta(M^+)$ is zero. It follows that the inclusion of  $Y^+$
into $Y$, induces an injection of $H_2(Y^+,\bq)$ into $H_2(Y,\bq)$.
Let $\em$ denote the matrix for the intersection pairing on $H_2(\cy)$ with
respect to the
generators for $L.$ It is not hard to see that $\em$ can be obtained from $M^+$
direct sum a diagonal matrix with $r+i$ minus two's and $\nu$ minus four's on
the diagonal
by a sequence of row and column operations of the form where one modifies one
row (or column) by a multiple of another. It follows that  $|\det (\em)| =
2^d|\det (M^+)|$. Thus  $|\det (2M^+)|= 2^{\rho^+-d}|\det (\em)|$.

Since $|\det (\em)|\ne 0,$ the inclusion of  $\ef$
into $\cy$ induces an injection of $H_2(\ef,\bq)$ into $ H_2(\cy,\bq)$.
It follows that $\dim L= \Bbb L.$ Thus $L$ is a lattice of maximal rank in
$H_+.$ Let $\det H_+$ denote the  absolute
value of the determinant of the matrix for the restriction of the nonsingular
intersection pairing on $H_2(\cy)$ to $H_+$ with respect to some basis for
$H_+$. We have that $|\det (\em)| =
b^2 \det (H_+)$ where $b$ is an integer.

 The proof of \cite {Kh2,Lemma 3.7} then shows
that $|\det(H_+)| = 2^{h^+}.$   The proof of \cite {Kh2,Lemma 3.7} refers to
\cite {Kh3,Lemma 2.4} which is not precisely stated  in the English
translation and is given without proof. This same lemma is stated precisely
and given with proof in \cite{W,Lemma 3.14}.
Thus $|\det (2M^+)|= 2^{h^+ +\r^+-d}b^2$.  \qed
\par

\enddocument